\documentclass[11pt]{article}
\usepackage{latexsym,amsmath,amstext,mathrsfs,amssymb,amsthm,amscd,amsopn,verbatim,graphics}

\usepackage{times}
\usepackage{indentfirst}
\numberwithin{equation}{section}

\newcommand{\ndash}{\nobreakdash-\hspace{0pt}}

\DeclareMathOperator{\maurer}{MC}

\DeclareMathOperator{\hkr}{HKR}
\DeclareMathOperator{\id}{id}

\DeclareMathOperator{\tange}{T}

\DeclareMathOperator{\GL}{GL}

\newcommand{\tspace}[1]{\tange\! #1}

\newcommand{\rd}{\overleftarrow{\partial}} 
\newcommand{\ld}{\overrightarrow{\partial}}

\newcommand{\ii}{{\mathrm{i}}}

\newcommand{\dd}{\mathrm{d}}

\DeclareMathOperator{\ad}{ad}

\DeclareMathOperator{\Hom}{Hom}

\newtheorem{Thm}{Theorem}[section]
\newtheorem{Prop}[Thm]{Proposition}
\newtheorem{Lem}[Thm]{Lemma}

\theoremstyle{remark}

\theoremstyle{definition}

\newtheorem{Def}[Thm]{Definition}
\newtheorem{Exa}[Thm]{Example}

\newtheorem*{Reg*}{Regularization procedure}

\newcommand{\braket}[2]{\left\langle{\,{#1}\,,\,{#2}\,}\right\rangle}

\newcommand{\Lie}[2]{{\left[{\,{#1}\,,\,{#2}\,}\right]}}
\newcommand{\SN}[2]{{\left[{\,{#1}\,,\,{#2}\,}\right]_{_{\sf SN}}}}
\newcommand{\Ger}[2]{{\left[{\,{#1}\,,\,{#2}\,}\right]_{_{\sf G}}}}

\newcommand{\starLie}[2]{{\left[{\,{#1}\,\overset\star,\,{#2}\,}\right]}}
\newcommand{\Poiss}[2]{\left\{{\,{#1}\,,\,{#2}\,}\right\}}

\newcommand{\de}{\partial}

\newcommand{\ee}{\epsilon}
\newcommand{\eee}{\varepsilon}
\newcommand{\s}[1]{{\sigma\hbox{\tiny{$(#1)$}}}}
\newcommand{\bbK}{{\mathbb{K}}}
\newcommand{\bbC}{{\mathbb{C}}}
\newcommand{\bbR}{{\mathbb{R}}}

\newcommand{\bbZ}{{\mathbb{Z}}}
\newcommand{\bbN}{{\mathbb{N}}}
\newcommand{\Lg}{\mathfrak{g}}
\newcommand{\Lh}{\mathfrak{h}}

\newcommand{\Lm}{\mathfrak{m}}

\newcommand{\calH}{\mathcal{H}}

\newcommand{\calC}{\mathcal{C}}

\newcommand{\calG}{\mathcal{G}}

\newcommand{\calD}{\mathcal{D}}
\newcommand{\calL}{\mathcal{L}}
\newcommand{\calI}{\mathcal{I}}
\newcommand{\calM}{\mathcal{M}}

\newcommand{\calV}{\mathcal{V}}

\newcommand{\tLg}{\Tilde{\Lg}}
\newcommand{\tQ}{\widetilde{Q}}
\newcommand{\tZ}{\widetilde{Z}}

\newcommand{\sfA}{{\mathsf{A}}}

\newcommand{\sfB}{{\mathsf{B}}}

\newcommand{\sfG}{{\mathsf{G}}}
\newcommand{\sfa}{{\mathsf{a}}}

\newcommand{\sfV}{{\mathsf{V}}}
\newcommand{\sfT}{{\mathsf{T}}}

\newcommand{\sfX}{{\mathsf{X}}}
\newcommand{\sfY}{{\mathsf{Y}}}

\newcommand{\sfZ}{{\mathsf{Z}}}

\newcommand{\dede}[1]{{\frac\de{\de{#1}}}}

\DeclareMathOperator{\sgn}{sgn}

\DeclareMathOperator{\Ker}{Ker}
\DeclareMathOperator{\Imm}{Im}

\newcommand{\eps}{{\vcenter{\hbox{\scriptsize{$[[\epsilon]]$}}}}}
\newcommand{\Cinf}[1]{{C^{\infty}\!\!\left(#1\right)}}
\newcommand{\Bskew}{{B_1^{^{_{-}} }}}
\newcommand{\miss}[2]{{\mathop{\hat #1}_{#2}}}
\newcommand{\tinyotimes}{{\vcenter{\hbox{\tiny{$\otimes$}}}}} 
\newcommand{\mt}[1]{{\hbox{\tiny{$#1$}}}}
\newcommand{\bisum}[2]{{\genfrac{}{}{0pt}{}{#1}{#2}}}
\newcommand{\iso}{{\stackrel{\sim}{\longrightarrow}}}
\newcommand{\LL}{$L_{\infty}$\ndash}
\newcommand{\ext}[2]{{\Lambda^{#2}\!\left(#1\right)}}
\newcommand{\sym}[2]{{S^{#2}\!\left(#1\right)}}
\newcommand{\Ext}[2]{{\overline{\Lambda}^{#2}\!\left(#1\right)}}
\newcommand{\Sym}[2]{{\overline{S}^{#2}\!\left(#1\right)}}
\newcommand{\Dbar}{{\overline{D}}}
\newcommand{\aff}[1]{{#1}^{\mathit{aff}}}
\newcommand{\coor}[1]{{#1}^{\mathit{coor}}}
\DeclareMathOperator{\Star}{Star}
\DeclareMathOperator{\Conf}{Conf}
\DeclareMathOperator{\dec}{dec}

\begin{document}

\title{{Formality and Star Products}}
\author{Alberto~S.~Cattaneo\\
Lecture notes taken by Davide Indelicato\\ \\
Institut f\"ur Mathematik,\\ Universit\"at Z\"urich--Irchel, \\
Winterthurerstrasse 190,\\ CH-8057 Z\"urich, Switzerland  \\ \\
\texttt{asc@math.unizh.ch}\\
\texttt{inde@math.unizh.ch}
}

\maketitle

\abstract{These notes, based on the mini-course
given at the PQR2003 Euroschool
held in Brussels in 2003,  
aim to review Kontsevich's formality
theorem
together with his
formula for the
star product on a given Poisson manifold.
A brief
introduction to the employed mathematical tools and physical motivations
is also given.}

\clearpage

\tableofcontents
\vfill
\section*{Acknowledgments}
We thank the hospitality and financial support
at Euroschool and Euroconference PQR2003.
A.~S.~C. acknowledges partial support of SNF Grant No.~20-100029/1.

We are very grateful to Carlo A. Rossi and 
Jim Stasheff for carefully reading 
the manuscript and making several suggestions to improve the text.
We also want to thank Benoit Dherin and Luca Stefanini for 
the stimulating discussions we have had with them.

\clearpage 

\section{Introduction} 

This work is based on the course given during the international Euroschool
on {\em Poisson Geometry, Deformation Quantisation
and Group Representations} held in Brussels in 2003.

The main goal is to describe Kontsevich's proof 
of the formality of the (differential graded) Lie
algebra of multidifferential operators on $\bbR^d$
and its relationship to the existence and classification
of star products on a given Poisson manifold. We
start with a survey of the physical background 
which gave origin to such a problem and a historical review
of the subsequent steps which led to the final solution.

\subsection*{Physical motivation} 

In this Section we give a brief overview of 
physical motivations that led to the genesis of
the deformation quantization problem, referring
to the next Sections and to the literature cited
throughout the paper for a precise definition of
the mathematical structures we introduce.

In the hamiltonian formalism of classical mechanics, 
a physical system is described by an even-dimensional
manifold $M$ --- the phase space --- endowed with 
a symplectic (or more generally a Poisson) structure
together with a smooth function $H$ ---
the hamiltonian function --- on it. A physical state
of the system is represented by a point in $M$ while the 
physical observables (energy, momentum and so on) correspond
to (real) smooth functions on $M$. The time evolution of
an observable $O$ is governed by an equation
of the form
\[
\frac{d\,O}{d\,t} = \Poiss{H}{O} 
\]
where $\Poiss{}{}$ is
the Poisson bracket on $\Cinf{M}$. This bracket
is completely determined by its action on
the coordinate functions
\[
\Poiss{p_i}{q_j} = \delta_{ij}
\]
(together with $ \Poiss{p_i}{p_j} = \Poiss{q_i}{q_j} = 0$)
where $(q_1, \ldots , q_n, p_1, \ldots , p_n)$ are 
local coordinates on the $2n$\ndash dimensional manifold $M$.

On the other hand, a quantum system is described
by a complex Hilbert space $\calH$ together with an operator $\widehat{H}$.
A physical state of the system is represented by a vector\footnote{
Actually, due to the linearity of the dynamical equations, 
there is a non-physical multiplicity which can be avoided
rephrasing the quantum formalism on a projective
Hilbert space, thus identifying a physical state with
a ray in $\calH$
} in $\calH$ while the physical observables are now 
self-adjoint operators in $\calL(\calH)$. 
The time evolution of such an operator in 
the Heisenberg picture is given by
\[
\frac{d\,\widehat{O}}{d\,t}  = \frac{i}{\hbar}\Lie{\widehat{H}}{\widehat{O}} 
\]
where $\Lie{}{}$ is the usual commutator which endows $\calL(\calH)$
with a Lie algebra structure. 
The correspondence with classical mechanics is completed
by the introduction of the position $\widehat{q}_i$ 
and momentum $\widehat{p}_j$ operators, which satisfy
the canonical commutation relations:
\[
\Lie{\widehat{p}_i}{\widehat{q}_j} = \frac{i}{\hbar} \delta_{ij}.
\]

This correspondence is by no means a mere analogy,
since quantum mechanic was born to replace the hamiltonian 
formalism in such a way that the classical picture could 
still be recovered as a ``particular case''.
This is a general principle in the development of a new 
physical theory: whenever experimental phenomena 
contradict an accepted theory, a new one is sought
which can account for the new data, but still reduces 
to the previous formalism when the new parameters 
introduced go to zero. In this sense, classical mechanics
can be regained from the quantum theory in the limit
where $\hbar$ goes to zero.

The following question naturally arises: is there
a precise mathematical formulation of this quantization 
procedure in the form of a well-defined map between classical objects
and their quantum counterpart?

Starting from the canonical quantization method for $\bbR^{2 n}$,
in which the central role is played by the
canonical commutation relation, a first approach was 
given by {\sf geometric quantization}:
the basic idea underlying this theory 
was to set a relation between the phase space $\bbR^{2n}$ 
and the corresponding Hilbert space $\calL(\bbR^{n})$
on which the Schr\"odinger equation is defined.
The  first works on geometric quantization are 
due to Souriau \cite{So}, Kostant \cite{Kos} and  Segal \cite{Se}, 
although many of their ideas were based on previous works by
Kirillov \cite{Kir}. We will not discuss further this approach,
referring the reader to the cited works.

On the other hand, one can focus attention on 
the observables instead of the physical states, 
looking for a procedure to get the non-commutative 
structure of the algebra of operators from the 
commutative one on $\Cinf{\bbR^{2n}}$. However, one of the
first result achieved was the ``no go'' theorem by Groenwold
\cite{Gro} which states the impossibility of quantizing
the Poisson algebra $\Cinf{\bbR^{2n}}$ in such a way 
that the Poisson bracket of any two functions 
is sent onto the Lie bracket of the two corresponding 
operators. Nevertheless, instead of mapping
functions to operators, one can ``deform''
the pointwise product on functions into a
non-commutative one, realizing, in an autonomous manner,
quantum mechanics directly on $\Cinf{\bbR^{2n}}$: this 
is the content of the {\sf deformation quantization} program promoted 
by Flato in collaboration with Bayen, Fronsdal, Lichnerowicz 
and Sternheimer,

\subsection*{Historical review of deformation quantization}

The origins of the deformation quantization approach
can be traced back to works of Weyl's \cite{We}, who gave 
an explicit formula for the operator $\Omega(f)$ on $\calL(\bbR^n)$
associated to a smooth function $f$ on the phase space
$\bbR^{2n}$:
\[
\Omega(f) := \int_{\bbR^{2n}} \check{f}(\xi, \eta)\;
e^{\frac{i}{\hbar}(P \cdot \xi + Q \cdot \eta)} d^n \xi\, d^n \eta ,
\]      
where $\check{f}$ is the inverse Fourier transform of $f$, 
$P_i$ and $Q_j$ are operators satisfying the canonical 
commutation relations and the integral is taken in the weak sense.
The arising problem of finding an inverse formula
was solved shortly afterwards by Wigner \cite{Wi}, who gave
a way to recover the classical observable from the quantum
one taking the symbol of the operator. It was then Moyal \cite{Mo}
who interpreted the symbol of the commutator of two 
operators corresponding to the functions $f$ and $g$
as what is now called a {\sf Moyal bracket} $\calM$:
\[
\calM(f, g) = \frac{\sinh (\ee\; P)}{\ee}(f, g) = 
\sum_{k=0}^{\infty} \frac{\ee^{2k} }{(2 k +1)!}\; P^{2 k+1}(f, g),
\] 
where $\ee = \frac{i \hbar}{2}$ and $P^k$ is the $k$\ndash th power
of the Poisson bracket on $\Cinf{R^{2n}}$. A similar formula
for the symbol of a product $\Omega(f) \Omega(g)$ had already been
found by Groenewold \cite{Gro} and can now be interpreted
as the first appearance of the Moyal star product $\star$, 
in terms of which the above bracket can be rewritten as
\[
\calM(f, g) = \frac{1}{2 \ee} (f \star g - g \star f).
\]

However, it was not until Flato gave birth to 
his program for deformation quantization that this 
star product was recognized as a non commutative deformation of
the (commutative) pointwise product on the algebra of
functions. This led to the first paper \cite{FLS1}
in which the problem was posed of giving a general recipe
to deform the product in $\Cinf{M}$ in such a way 
that $\frac{1}{2 \ee} (f \star g - g \star f)$
would still be a deformation of the given Poisson
structure on $M$. Shortly afterward Vey \cite{Ve}
extended the first approach, which considered only
$1$\ndash differentiable deformation, to more general 
differentiable deformations, rediscovering in an 
independent way the Moyal bracket. This opened the way 
to subsequent works (\cite{FLS2} and \cite{BFFLS}) in which 
quantum mechanics was formulated as a deformation
(in the sense of Gerstenhaber theory) of 
classical mechanics and the first significant applications
were found.

The first proof of the existence of star products on 
a generic symplectic manifold was given by DeWilde and 
Lecomte \cite{DL} and relies on the fact that locally 
any symplectic manifold of dimension $2 n$ 
can be identified with $\bbR^{2 n}$ via a Darboux chart.
A star product can thus be defined locally by the Moyal
formula and these local expressions can be glued 
together by using cohomological arguments.

A few years later and independently of this previous result, 
Fedosov \cite{Fed} gave an explicit algorithm to construct 
star products on a given symplectic manifold:
starting from a symplectic connection on $M$, he defined
a flat connection $D$ on the Weyl bundle associated
to the manifold, to which the local Moyal expression
for $\star$ is extended; the algebra of (formal) functions
on $M$ can then be identified with the subalgebra of horizontal
sections w.r.t.\ $D$.
We refer the reader to Fedosov's book for the details.
This provided a new proof of existence which could be 
extended to regular Poisson manifolds and opened the way to 
further developments.

Once the problem of existence was settled,
it was natural to focus on the classification
of equivalent star products, where the equivalence
of two star products has to be understood in the sense 
that they give rise to the same algebra up to 
the action of formal automorphisms which are deformations
of the identity. Several authors came
to the same classification result using very 
different approaches, confirming what was
already in the seminal paper \cite{BFFLS} by Flato et
al.\, namely that the obstruction to equivalence 
lies in the second de Rham cohomology of the manifold $M$.
For a comprehensive enumeration of the different proofs
we address the reader to \cite{DS}.

The ultimate generalization to the case of a generic
Poisson manifold relies on the formality theorem Kontsevich
announced in \cite{Ko1} and subsequently proved in \cite{Ko2}.
In this last work he derived  an explicit 
formula for a star product on $\bbR^d$, 
which can be used to define it locally on any $M$.
Finally, Cattaneo, Felder and Tomassini \cite{CFT1} gave a globalization 
procedure to realize explicitly what Kontsevich
proposed, thus completing the program 
outlined some thirty years before by Flato.

For a complete overview of the process which 
led from the origins of quantum mechanics to
this last result and over, we refer to 
the extensive review given by Dito and Sternheimer
in \cite{DS}.

As a concluding remark, we would like to mention
that the Kontsevich formula can also be expressed
as the perturbative expression of the functional integral of a
topological field theory --- the so-called Poisson sigma model 
(\cite{Ik}, \cite{SS}) --- 
as Cattaneo and Felder showed in \cite{CF1}. 
The diagrams Kontsevich introduced for his construction
of the local expression of the star product arise naturally
in this context as Feynman diagrams corresponding 
to the perturbative evaluation of a certain observable.

\subsection*{Plan of the work}

In the first Section we introduce the basic definition
and properties of the star product in the most general setting and
give the explicit expression of the Moyal product on $\bbR^{2 d}$
as an example. The equivalence relation on star products 
is also discussed, leading to the formulation of the classification
problem.

In the subsequent Section we establish the relation between 
the existence of a star product on a given manifold $M$ and the 
formality of the (differential graded) Lie algebra $\calD$
of multidifferential operators on $M$. We introduce the main tools
used in Kontsevich's construction and present the 
fundamental result of Hochschild, Kostant and Rosenberg 
on which the formality approach is based.

A brief digression follows, in which the 
formality condition is examined from a dual point 
of view. The equation that the formality map 
from the (differential graded) Lie algebra
$\calV$ of multivector fields to $\calD$ must fulfill
is rephrased in terms of an infinite family of equations 
on the Taylor coefficients of the dual map.

In the third Section Kontsevich's construction
is worked out explicitly and the formality 
theorem for $\bbR^d$ is proved following
the outline given in \cite{Ko2}. Finally,
the result is generalized to any Poisson manifold 
$M$ with the help of the globalization procedure
contained in \cite{CFT1}.

\section{The star product} \label{star-prod}

In this Section we will briefly give the definition and main 
properties of the star product.\\
Morally speaking, a star product is a formal non-commutative deformation of the usual 
pointwise product of functions on a given manifold.
To give a more general definition, one can start with a commutative associative 
algebra $\sfA$ with unity over a base ring $\bbK$ and deform it to the 
algebra $\sfA \eps$ over the ring of formal power series $\bbK\eps$.
Its elements are of the form
\[
\hspace{2.5cm}
C =  \sum_{i=0}^{\infty} c_i \, \ee^i 
\hspace{2cm}
c_i \in \sfA
\]
and the product is given by the Cauchy formula, multiplying
the coefficients according to the original product on $\sfA$
\[
\Big( \sum_{i=0}^{\infty} a_i \, \ee^i \Big) \bullet_{\ee}
\Big( \sum_{j=0}^{\infty} b_j \, \ee^j \Big)
=
\sum_{k=0}^{\infty} \Big( \sum_{l=0}^k a_{k-l} \cdot b_{l} \Big) \, \ee^k 
\]
The star product is then a  $\bbK\eps$\ndash linear associative product $\star$
on $\sfA \eps$ which deforms this trivial extension 
$\bullet_{\ee}\colon \sfA\eps \otimes_{\bbK\eps} \sfA\eps \to \sfA\eps$ in the
sense that for any two $v,w \in \sfA\eps$
\[
v \star w = v \bullet_{\ee} w \mod \ee.
\]

\par

In the following we will restrict our attention to the case in which $\sfA$ is 
the Poisson algebra $\Cinf{M}$ of smooth functions on $M$ endowed with the usual 
pointwise product
\[
\hspace{2.5cm}
f\cdot g(x) := f(x)\; g(x) 
\hspace{1.5cm}
\forall x \in M
\] 
and $\bbK$ is $\bbR$. 

With these premises we can give the following
\begin{Def} \label{def-star-prod}
A {\sf star product} on $M$ is an  $\bbR\eps$\ndash bilinear map
\[
\begin{array}{ccc}
\Cinf{M}\eps \times \Cinf{M}\eps &\to& \Cinf{M}\eps\\
(f, g) &\mapsto& f \star g
\end{array}
\]
such that

\begin{itemize}
\item[i)]   $ f \star g = f \cdot g +  \sum_{i=1}^{\infty} B_i(f,g)\;\ee^i$,
\item[ii)]  $ ( f \star g ) \star h = f \star (g \star h) 
\hspace{1.5cm}\forall f,g,h \in \Cinf{M}
\hspace{1cm}( associativity )$,
\item[iii)] $ 1 \star f = f \star 1 = f
\hspace{2.5cm}\forall f\in \Cinf{M}$.
\end{itemize}
The $B_i$ could in principle be just bilinear operators, 
but,  in order to encode locality from a physical point of view,
one requires them to be bidifferential operators on $\Cinf{M}$ 
of globally bounded order, that is, bilinear operators which 
moreover are differential operators w.r.t.\ each argument; 
writing the $i$\ndash th term in local coordinates: 
\[
B_i(f,g)= \sum_{K,L} \beta_i^{K L} \, \de_K f \, \de_L g 
\]
where the sum runs over all multi-indices $K = (k_1, \ldots, k_m)$
and $L=(l_1, \ldots, l_n)$ of any length $m,n \in \bbN$ and the usual
notation for higher order derivatives is applied; 
the $\beta_i^{K L}$'s are 
smooth functions, which are non-zero only for finitely many choices of the 
multi-indices $K$ and $L$.
\end{Def}

\begin{Exa} {\bf The Moyal star product}

We have already introduced the Moyal star product as the 
first example of a deformed product on the algebra of
functions on $\bbR^{2d}$ endowed with the canonical symplectic form.
Choosing Darboux coordinates $(q,p)=(q_1, \dots , q_d, p_1, \ldots , p_d)$
we can now give an explicit formula for the product of two functions 
$f,g \in \Cinf{\bbR^{2 d}}$:
\[
f \star g \;(q,p) 
:= f(q,p)\; \exp \left(\ii \frac{\hbar}{2}
\left( \rd_q \ld_p -  \rd_p \ld_q \right) \right)\; g(q,p),
\]
where the $\rd$'s operate on $f$ and the $\ld$'s on $g$; 
the parameter $\ee$ has been replaced by the expression
$\ii \frac{\hbar}{2}$ that usually appears in the physical literature.

More generally, given a constant skew-symmetric tensor
$\lbrace\alpha^{ij}\rbrace$ on $\bbR^d$ with $i,j=1, \ldots, d$, 
we can define a star product by:
\begin{equation} \label{gen-moyal}
f \star g\; (x) = \exp \left(\ii \frac{\hslash}{2}\; \alpha^{ij} \dede{x^i} \dede{y^j}  \right) 
f(x)\, g(y)\Big|_{y=x}.
\end{equation}

We can easily check that such a star product is associative for any choice of $\alpha_{ij}$
\[
\begin{aligned}
((f \star g) \star h) \; (x) 
&= e^{\left(\ii \frac{\hslash}{2}\; \alpha^{ij} \dede{x^i} \dede{z^j}  \right)} 
(f \star g) (x) h(z)  \Big|_{x=z}  =\\
&= e^{ \left(\ii \frac{\hslash}{2}\; \alpha^{ij} (\dede{x^i} + \dede{y^i}) \dede{z^j} \right)}
e^{\left(\ii \frac{\hslash}{2}\; \alpha^{kl} \dede{x^k} \dede{y^l}\right)} f(x) g(y) h(z) \Big|_{x=y=z} =\\
&= e^{ \left(\ii \frac{\hslash}{2}\; \alpha^{ij} \dede{x^i} \dede{z^j} +
\alpha^{kl} \dede{y^k} \dede{z^l}+\alpha^{mn} \dede{x^m} \dede{y^n}\right) } f(x) g(y) h(z) \Big|_{x=y=z} =\\
&= e^{ \left(\ii \frac{\hslash}{2}\; \alpha^{ij} \dede{x^i} (\dede{y^j} + \dede{z^j}) \right)}
e^{\left(\ii \frac{\hslash}{2}\; \alpha^{kl} \dede{y^k} \dede{z^l}\right)} f(x) g(y) h(z) \Big|_{x=y=z} =\\
&= (f \star (g \star h)) \; (x).
\end{aligned}
\]
Point $i)$ and $iii)$ in Definition \eqref{def-star-prod} and the $\bbR\eps$\ndash linearity
can be checked as well directly from the formula \eqref{gen-moyal}.
\end{Exa}

\bigskip\bigskip

We would like to emphasize that condition $iii)$ in the Definition \ref{def-star-prod}
implies that the degree 0 term in the r.h.s.\ of $i)$ has to be the usual product 
and it moreover ensures that the $B_i$'s are bidifferential operators
in the strict sense, i.e.\ they have no term of order $0$
\begin{equation} \label{Bi-on-const}
B_i(f,1) = B_i(1,f)=0\qquad\forall i\in \bbN_0 .
\end{equation} 

As another consequence of the previous requirements on the $B_i$'s, 
it is straightforward to prove that the skew-symmetric part $\Bskew$ of the 
first bidifferential operator, defined by
\[
\Bskew(f,g) := \frac{1}{2} \Big( B_1(f,g)- B_1(g,f) \Big)
\]
satisfies the following equations:
\begin{itemize}
\item[-] $\Bskew (f,g) = - \Bskew(g,f)$,
\item[-] $\Bskew(f,g \cdot h) = g \cdot \Bskew(f,h) + \Bskew(f,g)\cdot h$,
\item[-] $\Bskew (\Bskew(f,g),h) + \Bskew (\Bskew(g,h),f) + \Bskew (\Bskew(h,f),g) = 0$.
\end{itemize}
A bilinear operator on $\Cinf{M}$ which satisfies these three identities is called 
a {\sf Poisson bracket}. A smooth manifold $M$ endowed with 
a Poisson bracket on the algebra of smooth functions is called a 
{\sf Poisson manifold} (see also \cite{BW} and references therein).

It is therefore natural to look at the inverse problem: given a Poisson manifold $M$, 
can we define an associative, but possibly non commutative, product
$\star$ on the algebra of smooth functions, 
which is a deformation of the pointwise product and such that 
\[
\frac{f \star g - g \star f}{\ee} \mod \ee= \Poiss{f}{g}
\]
for any pair of functions $f,g \in \Cinf{M}$?

In order to reduce an irrelevant multiplicity of solutions,
the problem can be brought down to the study of 
equivalence classes of such products,
where the equivalence is to be understood in the sense of the following
\begin{Def} \label{star-prod-eq}
Two star products $\star$ and $\star'$ on $\Cinf{M}$
are said to be {\sf equivalent} if{f}  there exists a linear operator 
$\calD\colon \Cinf{M}\eps \to \Cinf{M}\eps$ of the form
\[
\calD f := f +   \sum_{i=1}^{\infty} D_i(f)\;\ee^i
\]
such that
\begin{equation} \label{eq-star-prod}
f \star' g = \calD^{-1} \left( \calD f \star \calD g \right) 
\end{equation}
where $\calD^{-1}$ has to be understood as the inverse in the sense of formal power series.
\end{Def}
It follows from the very definition of star product 
that also the $D_i$'s have to be differential operators
which vanish on constants, as was shown in \cite{GR}
(and without proof in \cite{Ve}).

This notion of equivalence leads immediately to a generalization 
of the previously stated problem, according to the following
\begin{Lem}
In any equivalence class of star products, there exists a representative
whose first term $B_1$ in the $\ee$ expansion is skew-symmetric.
\end{Lem}

\begin{proof}
Given any star product 
\[
f \star g := f \cdot g + \ee\; B_1(f,g) + \ee^2 B_2(f,g) + \cdots
\]
we can define an equivalent star product as in \eqref{eq-star-prod}
with the help of a formal differential operator
\[
\calD = \id + \ee\;  D_1 + \ee^2 D_2 + \cdots.
\]
The condition for the first term of the new star product to be skew-symmetric 
$B_1' (f,g) + B_1'(g,f)=0$
gives rise to an equation for the first term of the differential operator
\begin{equation} \label{d1}
D_1(f \; g) = D_1 f \; g + f \; D_1 g + 
\frac{1}{2} \Big( B_1(f,g) + B_1(g,f) \Big),
\end{equation}
which can be used to define $D_1$ locally on polynomials and hence by completion 
on any smooth function. By choosing a partition of unity, we may finally
apply $D_1$ to any smooth function on $M$.

We can start by choosing  $D_1$ to vanish on linear functions. 
Then the equation \eqref{d1} defines uniquely the action of $D_1$ 
on quadratic terms, given by the symmetric part $B_1^{+}$ of the bilinear operator $B_1$:
\[
D_1(x^i x^j) =B_1^{+} (x^i, x^j):= \frac{1}{2} \Big( B_1(x^i,x^j) + B_1(x^j,x^i) \Big).
\]
where $\lbrace x^k \rbrace$ are local coordinates on the manifold $M$.
The process extends to any monomial and ---
as a consequence of the associativity of $\star$ 
--- gives rise to a well defined operator
since it does not depend on the way we group the factors.
We check this on a cubic term:
\[
\begin{aligned}
D_1((x^i x^j)\, x^k) &= D_1(x^i x^j)\,  x^k +  x^i x^j\,  D_1(x_k) +B_1^{+} (x^i x^j, x^k) = \\
& = B_1^{+}( x^i, x^j) \, x^k + B_1^{+} (x^i x^j, x^k)= \\
& = B_1^{+}( x^i, x^j  x^k) + x^i \, B_1^{+} (x^j, x^k)= \\
& = x^i \, D_1(x^j x^k) + D_1(x_i)\, x^j x^k + B_1^{+} (x^i, x^j x^k) = D_1(x^i\, (x^j x^k)). 
\end{aligned}
\]
The equality between the second and the third lines is a consequence of
the associativity  of the star product: it is indeed the term of order $\ee$ in 
$(x^i \star x^j) \star x^k = x^i \star (x^j \star x^k) $
once we restrict the operators appearing on both sides to their symmetric part.
\end{proof}

The above proof is actually a particular case of the 
Hochschild--Kostant--Rosenberg theorem.
Associativity implies in fact that $B_1^{+}$ is a Hochschild cocycle, 
while in \eqref{d1} we want to express it as a Hochschild coboundary:
the HKR theorem states exactly that this is always possible
on $\bbR^d$ and thus locally on any manifold.

From this point of view, the natural subsequent step is
to look for the existence and uniqueness of equivalence classes of star products
which are deformations of a given Poisson structure on the smooth manifold $M$. 
As already mentioned in the introduction, the existence of such products was 
first proved by DeWilde and Lecomte \cite{DL} in the symplectic case,
where the Poisson structure is defined via a symplectic form 
(a non degenerate closed $2$\ndash form).
Independently of this previous result, Fedosov \cite{Fed}
gave an explicit geometric construction: the star product 
is obtained  ``glueing'' together local expressions obtained via the Moyal formula.

As for the classification,
the role played by the second de Rham cohomology of the manifold, 
whose occurrence in connection with this problem can be traced back to 
\cite{BFFLS}, has been clarified in subsequent works by different authors 
(\cite{NT}, \cite{BCG}, \cite{Gu}, \cite{Xu}, \cite{Bon}, \cite{De})
until it came out that equivalence classes of star products 
on a symplectic manifold are in one-to-one correspondence with elements
in $H^2_{dR}(M)\eps$.

The general case was solved by Kontsevich in \cite{Ko2}, who gave an explicit
recipe for the construction of a star product starting from any
Poisson structure on $\bbR^d$. This formula can thus be used to define
locally a star product on any Poisson manifold; the local expressions
can be once again glued together to obtain a global star product, as 
explained in Section \ref{globalization}.
As already mentioned, this result
is a straightforward consequence of the formality theorem, which 
was already announced as a conjecture in 
\cite{Ko1} and subsequently proved in \cite{Ko2}.
In the following, we will review this stronger result which relates two 
apparently very different mathematical objects --- multivector fields 
and multidifferential operators --- and we will come to the explicit 
formula as a consequence in the end.

\bigskip

As a concluding act,  we anticipate the Kontsevich formula
even though we will fully understand its meaning only in the forthcoming Sections.   

\begin{equation}
f \star g := f \cdot g +
\sum_{n=1}^{\infty} \ee^n \hspace{-.3cm} \sum_{\Gamma \in G_{n,2}} w_{\Gamma}\; B_{\Gamma} (f,g)
\end{equation}
The bidifferential operators as well as the weight coefficients 
are indexed by the elements $\Gamma$ of a suitable subset $G_{n,2}$
of the set of graphs on $n+2$ vertices, the so-called 
{\sf admissible graphs}.

\section{Rephrasing the main  problem: the formality}

In this Section we introduce the main tools that we will need to review 
Kontsevich's construction of a star product on a Poisson manifold.

The problem of classifying star products on a given Poisson 
manifold $M$  is solved by proving that there is a one-to-one 
correspondence between equivalence classes of star products and equivalence classes
of formal Poisson structures.

While the former were defined in the previous Section,
the equivalence relation on the set of formal Poisson structures is defined as follows.
First of all, to give a Poisson structure on $M$ is the same as to choose 
a {\sf Poisson bivector field},
i.e.\ a section $\pi$ of  $\bigwedge^2 \tspace{M}$ with certain properties 
that we will specify later, and define the Poisson bracket via the pairing between
(exterior powers of the) tangent and cotangent space:
\begin{equation} \label{p-bivector}
\hspace{3.5 cm}
\Poiss{f}{g}:= \frac{1}{2}\;\braket{\pi}{\dd f \wedge \dd g} 
\qquad \forall f,g \in \Cinf{M}.
\end{equation}

The set of Poisson structures is acted on by the group of diffeomorphisms of $M$,
the action being given through the push-forward by
\begin{equation} \label{phi-on-pi}
\pi_{\phi} := \phi_{*} \pi.
\end{equation}

To extend this notion to formal power series, we can introduce a bracket
on $\Cinf{M}\eps$ by:
\begin{equation} \label{formal-poisson}
\Poiss{f}{g}_{\ee}:= \sum_{m=0}^{\infty} \ee^m \hspace{-0.4cm} 
\mathop{\sum_{i,j,k=0}^{m}}_{i+j+k=m}\braket{\pi_i}{\dd f_j \wedge \dd g_k}
\end{equation}
where 
\[
f = \sum_{j=0}^{\infty} \ee^j f_j
\qquad\text{and}\qquad 
g = \sum_{k=0}^{\infty} \ee^k g_k
\]
One says that
\[
\pi_{\ee}:= \pi_0 + \pi_1 \, \ee + \pi_2 \, \ee^2 + \cdots
\]
is a {\sf formal Poisson structure} if $\Poiss{}{}_{\ee}$ is a Lie bracket 
on $\Cinf{M}\eps$.

The gauge group in this case is given by {\sf formal diffeomorphisms},
i.e.\ formal power series of the form
\[
\phi_{\ee}:= \exp(\ee \, \sfX)
\]
where $\sfX:=\sum_{k=0}^{\infty} \ee^k X_k $ is a {\sf formal vector field}, i.e.\
a formal power series whose coefficients are vector fields.
This set is given the structure of a group defining the product 
of two such exponentials via the Baker--Campbell--Hausdorff formula:
\begin{equation} \label{bch}
\exp(\ee \, \sfX) \cdot  \exp(\ee \, \sfY) :=
\exp(\ee\;\sfX + \ee\;\sfY + \frac{1}{2} \ee\;\Lie{\sfX}{\sfY} + \cdots).
\end{equation}
The action which generalizes \eqref{phi-on-pi} is then given
via the Lie derivatives $\calL$ on bivector fields by
\begin{equation} \label{gauge-action}
\exp \left( \ee\, \sfX \right)_{*} \pi := 
\sum_{m=0}^{\infty} \ee^m \hspace{-0.4cm} 
\mathop{\sum_{i,j,k=0}^{m}}_{i+j+k=m} (\calL_{\sfX_i})^j \pi_k
\end{equation}

\bigskip

Kontsevich's main result in \cite{Ko2} was to find an identification between the
set of star products modulo the action of the differential operators defined in 
\eqref{star-prod-eq} and the set of formal Poisson structure modulo this gauge group.
(For further details the reader is referred to \cite{Arb} and \cite{Ma})

\subsection{DGLA's, L$_{\infty}$- algebras and deformation functors} \label{dgla}

In the classical approach to deformation theory, (see e.g \cite{Art})
to each deformation is attached a DGLA via the solutions to  the Maurer--Cartan 
equation modulo the action of a gauge group.
The first tools we need to approach our problem are then contained in the following
definitions.
\begin{Def} \label{def-gla}
A {\sf graded Lie algebra} (briefly GLA) is a $\bbZ$\ndash graded
vector space $\Lg = \bigoplus_{i \in \bbZ} \Lg^i$
endowed with a bilinear operation
\[
\Lie{}{} \colon \Lg \otimes \Lg \to \Lg
\]
satisfying the following conditions:
\begin{itemize}
\item[a)] $\Lie{a}{b} \in \Lg^{\alpha+\beta}$ 
\hspace{\stretch{1}}
({\em homogeneity})
\item[b)] $\Lie{a}{b} = - (-)^{\alpha \beta} \Lie{b}{a}$ 
\hspace{\stretch{1}}
({\em skew-symmetry})
\item[c)] $\Lie{a}{\Lie{b}{c}} = \Lie{\Lie{a}{b}}{c} +(-)^{\alpha \beta} \Lie{b}{\Lie{a}{c}}$ 
\hspace{\stretch{1}}
({\em Jacobi identity})
\end{itemize}
for any $a \in \Lg^{\alpha}$, $b \in \Lg^{\beta}$ and $c \in \Lg^{\gamma}$
\end{Def}

As an example we can consider any Lie algebra as a GLA concentrated in degree
0. Conversely, for any GLA $\Lg$, its degree zero part $\Lg^0$ 
(as well as the even part $\Lg^{even} := \bigoplus_{i \in \bbZ} \Lg^{2 i}$) is a Lie 
algebra in the usual sense.

\begin{Def} \label{def-dgla}
A {\sf differential graded Lie algebra} is a GLA $\Lg$ together with a differential, 
$\dd \colon \Lg \to \Lg$, i.e.\ a linear operator of degree 1 
($ \dd \colon \Lg^i \to \Lg^{i+1}$) which satisfies the Leibniz rule
\[
\qquad\qquad\qquad\qquad
\dd \Lie{a}{b}  = \Lie{\dd\,a}{b} + (-)^{\alpha} \Lie{a}{\dd\,b}
\qquad\qquad a \in \Lg^{\alpha}, b \in \Lg^{\beta}
\]
and squares to zero ($\dd \circ \dd =0$).
\end{Def}
Again we can make any Lie algebra into a DGLA concentrated in degree $0$ with trivial 
differential $\dd = 0$. More examples can be found for instance in \cite{Ma}. 
In the next Section we will introduce the two DGLA's that 
play a role in deformation quantization.

The categories of graded and differential graded Lie algebras are completed 
with the natural notions of morphisms as graded linear maps
which moreover commute with the differentials and the brackets
\footnote{
We recall that a graded linear map $\phi\colon \Lg \to \Lh$ 
of degree $k$ is a linear map such that 
$\phi(\Lg^i) \subset \Lh^{i+k}$ $\forall i \in \bbN$.
We remark that, in the case of DGLA's, a morphism has to be a degree $0$
linear map in order to commute with the other structures.
}.
Since we have a differential, we can form a cohomology complex out of any DGLA 
defining the cohomology of $\Lg$ as 
\[ 
\calH^i(\Lg) := \Ker(\dd \colon \Lg^i \to \Lg^{i+1}) \Big/ \Imm(\dd \colon \Lg^{i-1} \to \Lg^i).
\]
The set $\calH := \bigoplus_i \calH^i(\Lg)$ has a natural structure of graded vector 
space and, because of the compatibility condition between the differential $\dd$ 
and the bracket on $\Lg$, it inherits the structure of a GLA, defined 
unambiguously on equivalence classes $|a|, |b| \in \calH$ by:
\[
\Lie{|a|}{|b|}_{\calH} := \left| \Lie{a}{b}_{\Lg} \right|.
\]
Finally, the cohomology of a DGLA  can itself be turned into a DGLA
with zero differential.

It is evident that every morphism $\phi\colon \Lg_1 \to \Lg_2$ of DGLA's 
induces a morphism $\calH(\phi) \colon \calH_1 \to \calH_2$ between
cohomologies. Among these, we are particularly interested in the so-called
{\sf quasi-isomorphisms}, i.e.\ morphisms of DGLA's inducing isomorphisms 
in cohomology. Such maps generate an equivalence relation:
two DGLA's $\Lg_1$ and $\Lg_2$ are called {\sf quasi-isomorphic}
if they are equivalent under this relation.\footnote{
We want to stress the fact that the existence of a quasi-isomorphism
$\phi\colon \Lg_1 \to \Lg_2$ does not imply the existence of a 
``quasi-inverse''  $\phi^{-1} \colon \Lg_2 \to \Lg_1$: therefore
these maps do not define automatically an equivalence relation.
This is the main reason why we have to consider the broader category 
of \LL algebras.
}

\begin{Def} \label{def-formal}
A differential graded Lie algebra $\Lg$ is called {\sf formal} if 
it is quasi-isomorphic to its cohomology, regarded
as a DGLA with zero differential and the induced bracket.
\end{Def}

The main result of Kontsevich's work --- the {\sf formality theorem} 
contained in \cite{Ko2} -- was to show that the DGLA of multidifferential 
operators, which we are going to introduce in the next Section, is formal. 

In order to achieve this goal, however, one has to rephrase 
the problem in a broader category, which we will define in this Section, 
though its structure will become clearer in Section \ref{digression},  
where it will be analyzed from a dual point of view.

To introduce the notation that will be useful throughout,
we start from the very basic definitions.
\begin{Def}
A {\sf graded coalgebra} (briefly GCA in the following) on the base ring $\bbK$
is a $\bbZ$\ndash graded vector space $\Lh = \bigoplus_{i \in \bbZ} \Lh^i$
endowed with a comultiplication, i.e.\ a graded linear map  
\[
\Delta \colon \Lh \to \Lh \otimes \Lh
\]
such that 
\[
\Delta(\Lh^i) \subset \bigoplus_{j+k=i} \Lh^j \otimes \Lh^k
\]
and which moreover satisfies the coassociativity condition
\[
(\Delta \otimes \id) \Delta (a) = (\id \otimes \Delta) \Delta (a)
\]
for every $a \in \Lh$.
It is said to be {\sf with counit} if there exists a morphism
\[
\ee \colon \Lh \to \bbK
\]
such that $\ee (\Lh^i) = 0 $ for any $i>0$ and
\[
(\ee \otimes \id) \Delta (a) = (\id \otimes \ee) \Delta (a) = a
\] 
for every $a \in \Lh$.
It is said to be {\sf cocommutative} if 
\[
\sfT \circ \Delta = \Delta
\]
where $\sfT\colon \Lh \otimes \Lh \to \Lh \otimes \Lh$ is the twisting map, 
defined on a product $x \otimes y$ of homogeneous elements of degree 
respectively $|x|$ and $|y|$ by
\[
\sfT ( x \otimes y ) := (-)^{|x| |y|}\, y \otimes x
\]
and extended by linearity.
\end{Def}

Given a (graded) vector space $V$ over $\bbK$, we can define new graded vector 
spaces over the same ground field by:
\begin{equation} \label{tensor-alg}
\begin{array}{lr}
        \begin{array}{rcl}
        T(V) &:=& \bigoplus_{n=0}^{\infty} V^{\otimes n} \\
        \,&\,&\\
        \overline{T}(V) &:=& \bigoplus_{n=1}^{\infty} V^{\otimes n} 
        \end{array}
& \qquad
V^{\otimes n}:= \left\lbrace\begin{array}{lr}
\underbrace{V \otimes \cdots \otimes V}_{n} & n\geq1\\
\bbK & n=0
\end{array}\right.,
\end{array}
\end{equation}
and turn them into associative algebras
w.r.t.\ the tensor product. $T(V)$ has also a unit given by $1 \in \bbK$.
They are called respectively the {\sf tensor algebra} and the {\sf reduced tensor algebra}.
As a graded vector space, $T(V)$ can be endowed with a coalgebra structure defining the 
comultiplication $\Delta_{T}$ on homogeneous elements by:
\[
\begin{aligned}
\Delta_{T} (v_1 \tinyotimes \cdots \tinyotimes v_n) := 
&\; 1 \otimes(v_1 \tinyotimes \cdots \tinyotimes v_n)\\
+& \;\sum_{j=1}^{j=n-1} (v_1 \tinyotimes  \cdots \tinyotimes v_j)
\otimes (v_{j+1} \tinyotimes \cdots \tinyotimes  v_n)\\
+& \;(v_1 \tinyotimes \cdots \tinyotimes v_n) \otimes 1
\end{aligned}
\]
and the counit $\ee_{T}$ as the canonical projection 
$\ee_{T} \colon T(V) \to V^{\otimes 0} = \bbK$.
The projection $T(V) \stackrel{\overline{\pi}}{\rightarrow} \overline{T}(V)$ 
and the inclusion $\overline{T}(V) \stackrel{i}{\hookrightarrow} T(V)$ 
induce a comultiplication also on the reduced algebra,
which gives rise to a coalgebra without counit. 

The tensor algebra gives rise to two other special algebras, 
the {\sf symmetric} $\sym{V}{}$ and {\sf exterior} $\ext{V}{}$ algebras,
defined as vector spaces as the quotients of $T(V)$ 
by the two-sided ideals --- respectively $\calI_{S}$ and $\calI_{\Lambda}$ ---
generated by homogeneous elements of the form $v\otimes w - \sfT(v\otimes w)$
and $v\otimes w + \sfT(v\otimes w)$. These graded vector spaces
inherit the structure of associative algebras w.r.t.\ 
the tensor product.
The reduced versions $\Sym{V}{}$
and $\Ext{V}{}$ are defined replacing $T(V)$ by the
reduced algebra $\overline{T}(V)$.

Also in this case, the underlying vector spaces can be endowed 
with a comultiplication which gives them 
the structure of coalgebras (without counit in the reduced cases). 
In particular on  $S(V)$  the comultiplication is 
given on homogeneous elements $v \in V$ by 
\[
\Delta_{S} (v) := 1 \otimes v + v \otimes 1,
\]
and extended as an algebra homomorphism w.r.t.\ the tensor product.

All the usual additional structures that can be put on an algebra 
can be dualized to give a dual version on coalgebras.
Having in mind the structure of DGLA's, we introduce the analog
of a differential by defining first coderivations.

\begin{Def}
A {\sf coderivation} of degree $k$ on a GCA $\Lh$ is a graded linear map 
$\delta\colon \Lh^{i} \to \Lh^{i+k}$ which satisfies the (co--)Leibniz 
identity:
\[
\Delta \delta (v) = ( \delta \otimes \id) \Delta (v)
 + ((-)^{k |v| } \id \otimes \delta) \Delta (v)
\qquad \forall v \in \Lh^{|v|}
\]
A {\sf differential} $Q$ on a coalgebra is a coderivation of degree one
that squares to zero.
\end{Def}

With these premises, we can give the definition of the main object we will 
deal with.

\begin{Def}

An {\sf \LL algebra} is a graded vector space $\Lg$ on $\bbK$ endowed
with a degree 1 coalgebra differential $Q$ on the reduced symmetric space 
$\overline{S}(\Lg[1])$.\footnote{
We recall that given any graded vector space $\Lg$, we can
obtain a new graded vector space $\Lg[k]$ by shifting each component
by $k$, i.e.\ 
\[
\Lg[k] = \bigoplus_{i \in \bbZ} \Lg[k]^i
\qquad \text{where} \qquad 
\Lg[k]^i:= \Lg^{i+k}.
\]
}An  {\sf \LL morphism} 
$F \colon (\Lg, Q) \to (\tLg, \tQ)$ is a morphism 
\[
F \colon \Sym{\Lg[1]}{} \longrightarrow \Sym{\tLg[1]}{} 
\]
of graded coalgebras (sometimes called pre\ndash\LL morphism),
which moreover commutes with the differentials ($F Q = \tQ F$).
\end{Def}
As in the dual case an algebra morphism $f\colon \sym{\sfA}{} \to \sym{\sfA}{}$ 
(resp. a derivation $\delta \colon \sym{\sfA}{} \to \sym{\sfA}{}$) is uniquely determined 
by its restriction to an algebra $\sfA=\sym{\sfA}{1}$ because of the
homomorphism condition $f(a b) = f(a) f(b)$ (resp. the Leibniz rule),
an \LL morphism $F$ and a coderivation $Q$ are uniquely
determined by their projection onto the first component $F^1$ resp. $Q^1$.
It is useful to generalize this notation introducing the symbol
$F^i_j$ (resp. $Q^i_j$) for the projection to the $i$\ndash th
component of the target vector space restricted to the $j$\ndash th 
component of the domain space.\footnote{
With the help of this decomposition,
it can be showed that for any given $j$,
only finitely many $F^i_j$ (and analogously $Q^i_j$)
are non trivial, namely $F^i_j =0$ for $i>j$. For an explicit
formula we refer the reader to \cite{Gra} and \cite{C}.
} With this notation, we can express 
in a more explicit way the condition which $F$ 
(resp. $Q$) has to satisfy to be an \LL morphism 
(resp. a differential). Since, with the above notation, 
$Q Q$, $F Q$ and $\tQ F$ are coderivations (as it can be checked
by a straightforward computation), it is sufficient to verify these
conditions on their projection to the first component. 

We deduce that a coderivation $Q$ is a differential iff 
\begin{equation} \label{QQ}
\sum_{i=1}^{n} Q^1_i Q^i_n = 0 \qquad \forall n \in \bbN_0
\end{equation}
while a morphism $F$ of graded coalgebras is an  \LL morphism
iff
\begin{equation} \label{FQ-QF}
\sum_{i=1}^{n} F^1_i Q^i_n \;=\; \sum_{i=1}^{n} \tQ^1_i F^i_n
\qquad \forall n \in \bbN_0.
\end{equation}

In particular, for $n=1$ we have 
\[
Q^1_1 Q^1_1 = 0 
\qquad \text{and} \qquad 
F^1_1 Q^1_1 \;=\; \tQ^1_1 F^1_1;
\]
therefore every coderivation $Q$ induces the structure of a complex of 
vector spaces on $\Lg$ and every \LL morphism
restricts to a morphism of complexes $F^1_1$. 
We can thus generalize the definitions given for a DGLA to this case,
defining a quasi-isomorphism of \LL algebras to be an
\LL morphism $F$ such that $F^1_1$ is a
quasi-isomorphism of complexes. The notion of formality can be 
extended in a similar way. We quote a result on 
\LL quasi-isomorphisms we will need later,
which follows from a classification theorem on \LL algebras.
\begin{Lem} \label{quasi-inv}
Let $F \colon (\Lg, Q) \to (\tLg, \tQ)$ be an \LL morphism.
If $F$ is a quasi-isomorphism it admits a {\sf quasi-inverse},
i.e.\ there exists an \LL morphism  $G \colon (\tLg, \tQ) \to (\Lg, Q)$ 
which induces the inverse isomorphism in the corresponding 
cohomologies.
\end{Lem}

For a complete proof of this Lemma together with an explicit
expression of the quasi-inverse  and a discussion of the above
mentioned classification theorem we refer the reader to \cite{C}.\\

In particular, Lemma \ref{quasi-inv} implies that  
\LL quasi-isomorphisms define equivalence relations,
i.e.\ two \LL algebras are \LL quasi-isomorphic iff
there is an \LL quasi-isomorphism between them. This is 
considerably simpler then in the case of DGLA's, where
the equivalence relation is only generated by the 
corresponding quasi-isomorphisms, and explains finally why
\LL algebras are a preferred tool in the solution
of the problem at hand.

\begin{Exa} \label{DGLA-LL}
To clarify in what sense we previously introduced \LL algebras as 
a generalization of DGLA's, we will show how to induce an \LL algebra
structure on any given DGLA $\Lg$.

We have already a suitable candidate for $Q^1_1$, since we know that 
it fulfills the same equation as the differential $\dd$: we may
then define $Q^1_1$ to be a multiple of the differential.
If we write down explicitly \eqref{QQ} for $n=2$, we get:
\[
Q^1_1 \, Q^1_2 + Q^1_2 \, Q^2_2 = 0;
\]
since every $Q^i_j$ can be expressed in term of a combination of products
of some $Q^1_k$, $Q^2_2$ must be a combination of $Q^1_1$ acting on the first
or on the second argument of $Q^1_2$ (for an explicit expression of the general 
case see \cite{Gra}). Identifying $Q^1_1$ with $\dd$ (up to a sign),
the above equation has thus the same form as the 
compatibility condition between the bracket $\Lie{}{}$ and the differential
and suggests that $Q^1_2$ should be defined in terms of the Lie bracket.
A simple computation points out the right signs, 
so that the coderivation is completely determined by
\[
\begin{array}{lll}
Q^1_1 (a) := (-)^{\alpha} \dd a 
&\,&  a \in \Lg^{\alpha},\\
Q^1_2(b\, c) := (-)^{\beta(\gamma -1)} \Lie{b}{c}
&\,&
b \in \Lg^{\beta}, c \in \Lg^{\gamma},\\
Q^1_n = 0 
&\,& \forall n\geq3. 
\end{array}
 \]
The only other equation involving non trivial terms follows from \eqref{QQ}
when $n=3$:
\[
Q^1_1 \, Q^1_3 + Q^1_2 \, Q^2_3 + Q^1_3 \, Q^3_3 =0. 
\]
Inserting the previous definition and expanding $Q^2_3$ in terms of $Q^1_2$ we get
\begin{equation} \label{Q-jacobi}
\begin{aligned}
(-)^{(\alpha + \beta)(\gamma-1)} & 
\Lie{(-)^{\alpha(\beta-1)} \Lie{a}{b}}{c} +\\
(-)^{(\alpha + \gamma)(\beta-1)}(-)^{(\gamma-1)(\beta-1)}  & 
\Lie{ (-)^{\alpha(\gamma-1)} \Lie{a}{c}}{b} + \\
(-)^{(\beta + \gamma)(\alpha-1)} (-)^{(\beta + \gamma)(\alpha-1)} & 
\Lie{(-)^{\beta (\gamma-1)} \Lie{b}{c}}{a} = 0,
\end{aligned}
\end{equation}
which, after a rearrangement of the signs, turns out to be the (graded) Jacobi identity.

According to the same philosophy, a DGLA morphism $F \colon \Lg \to \tLg$
induces an \LL morphism $\overline{F}$ which is completely 
determined by its first component $\overline{F}_1^1 := F$. 
In fact, the only two non trivial conditions on $\overline{F}$ 
coming from \eqref{FQ-QF} with $n=0$ resp. $n=1$ are:
\[
\begin{aligned}
\overline{F}_1^1 Q_1^1(f) = \tQ_1^1 \overline{F}_1^1 (f)
&\Leftrightarrow 
F (d\, f) = \tilde{d}\, F(f) \\ 
\overline{F}_1^1 Q^1_2 (f g) + \overline{F}_1^2 Q^2_2 (f g) = 
\tQ^1_1 \overline{F}_1^2 (f g) + \tQ^1_2 \overline{F}_2^2  (f g) 
&\Leftrightarrow
F \Big( \Lie{f}{g} \Big) = \Lie{F(f)}{F(g)}
\end{aligned}
\]

\end{Exa}

If we had chosen $Q^1_3$ not to vanish, the identity \eqref{Q-jacobi} 
would have been fulfilled up to homotopy, i.e.\ up to a term of the form 
\[
\dd \rho(g,h,k) \pm \rho(\dd g,h,k) \pm \rho(g,\dd h,k) \pm \rho(g,h,\dd k),
\]
where $\rho \colon \Lambda^3 \Lg \to \Lg[-1]$; 
in this case $\Lg$ is said to have the 
structure of a {\sf homotopy Lie algebra}. 

This construction can be generalized, introducing the
canonical isomorphism between the symmetric and exterior algebra
(usually called {\sf d\'ecalage isomorphism}\footnote{
More precisely, the d\'ecalage isomorphism is
given on the $n$\ndash symmetric power of $\Lg$
shifted by one by
\[
\begin{aligned}
\dec_n \colon \sym{\Lg[1]}{n} &\to \ext{\Lg}{n}[n]\\
x_1 \cdots x_n & \mapsto (-1)^{\sum_{i=1}^n (n-i)(|x_i|-1)} x_1 \wedge \ldots \wedge x_n,
\end{aligned}
\]
where the sign is chosen precisely to compensate for
the graded antisymmetry of the wedge product.
}) to define for each $n$
a {\sf multibracket} of degree $2 -n$
\[
[\cdot  , \cdots , \cdot]_n \colon \Lambda^n \Lg \to \Lg[2-n]
\]
starting from the corresponding $Q^1_n$. Equation \eqref{QQ} gives rise
to an infinite family of condition on these multibracket. A graded
vector space $\Lg$ together with such a family of operators is 
a {\sf strong homotopy Lie algebra}(SHLA).

To conclude this overview of the main tools we will need in the following
--- and to give an account of the last term in the title of this Section ---
we introduce now the Maurer--Cartan equation of a DGLA $\Lg$:
\begin{equation} \label{MC-eq}
\dd\, a + \frac{1}{2} \Lie{a}{a} = 0 \qquad a \in \Lg^1,
\end{equation}
which plays a central role in deformation theory, as will
exemplified in next Section, in \eqref{MC-eq-calV} and \eqref{MC-eq-calD}.

It is a straightforward application of the definition \ref{def-gla}
to show that the set of solutions to this equation is preserved under 
the action of any morphism of DGLA's and --- as we will see in the next
Section --- of any \LL morphism between the corresponding \LL algebras.

There is another group which preserve the solutions to the Maurer--Cartan 
equation, namely the gauge group that can be defined canonically starting 
from the degree zero part of any formal DGLA.

It is a basic result of Lie algebra theory that there exists a 
functor $\exp$ from the category of nilpotent Lie algebras to the
category of groups. For every such Lie algebra 
$\Lg$, the set defined formally as $\exp(\Lg)$ 
can be endowed with the structure of a group defining the product
via the Baker--Campbell--Hausdorff formula as in \eqref{bch};
the definition is well-posed since the nilpotency ensures that 
the infinite sum reduces to a finite one.

In the case at hand, generalizing what was somehow anticipated 
in \eqref{bch}, we can introduce the formal counterpart $\Lg\eps$
of any DGLA $\Lg$ defined as a vector space by 
$\Lg\eps := \Lg \otimes \bbK\eps$ and show that it has the natural
structure of a DGLA.
It is clear that the degree zero part $\Lg^0\eps$ is a Lie algebra, 
although non--nilpotent. Nevertheless, 
we can define the gauge group formally as the set 
$\sfG := \exp( \ee\, \Lg^0\eps)$ and introduce a well--defined 
product taking the Baker--Campbell--Hausdorff formula
as the definition of a formal power series.  
Finally, the action of the group on $\ee\, \Lg^1\eps$ can 
be defined generalizing the adjoint action in \eqref{bch}.

Namely:
\[
\begin{aligned}
\exp(\ee\, g) \sfa :=& \sum_{n=0}^{\infty} \frac{(\ad g)^n}{n !} (\sfa) -
\sum_{n=0}^{\infty} \frac{(\ad g)^n}{(n+1) !}(\dd g) \\
=&\; \sfa + \ee\, \Lie{g}{\sfa} - \ee\; \dd g + o(\ee^2)
\end{aligned}
\]
for any $g \in \Lg^0\eps$ and $\sfa \in \Lg^1\eps$.

It is a straightforward computation to show that this action preserves
the subset $\maurer(\Lg) \subset \ee\, \Lg^1\eps$ of solutions to the 
(formal) Maurer--Cartan equation.\\

\subsection{Multivector fields and multidifferential operators}

As we already mentioned, a Poisson structure is completely defined 
by the choice of a bivector field satisfying certain properties; on the other hand 
a star product is specified by a family of bidifferential operators. 
In order to work out the correspondence between these two objects, we are 
finally going to introduce the two DGLA's they belong to: multivector 
fields $\calV$ and multidifferential operators $\calD$.

\subsubsection{The DGLA $\calV$} \label{dgla-calV}

A $k$\ndash multivector field $\sfX$ is a Section of the $k$\ndash th exterior
power  $\bigwedge^k \tspace{M}$ of the tangent space $\tspace{M}$; choosing 
local coordinates $\lbrace x^i \rbrace_{_{ i=1, \ldots, \dim M}}$ and 
denoting by $\lbrace \de_i \rbrace_{_{i=1, \dots, \dim M}}$
the corresponding basis of the tangent space:
\[
\sfX = \sum_{i_1,\ldots,i_k = 1}^{\dim M} X^{ i_1 \cdots i_k} (x) \;
\de_{i_1} \wedge \cdots \wedge \de_{i_k}.
\]
The direct sum of such vector spaces has 
thus the natural structure of a graded vector space
\[
\hspace{2 cm}
\widetilde{\calV} := \bigoplus_{i=0}^{\infty} \widetilde{\calV}^{i} 
\hspace{2 cm}
\widetilde{\calV}^{i} :=
\left\lbrace\begin{array}{lr}
\Cinf{M} & \qquad i = 0 \\
\Gamma(\bigwedge^i \tspace{M}) & i \geq 1
\end{array}\right. ,
\]
having added smooth functions in degree 0.

The most natural way to define a Lie structure on $\widetilde\calV$ 
is by extending the usual Lie bracket on vector fields given 
in terms of the Lie derivative w.r.t.\ the first vector field:
\[
\Lie{\sfX}{\sfY} := \calL_{\sfX} \sfY.
\]
The same definition can be applied to the case when the second argument
is a function, setting:
\[
\Lie{\sfX}{f} := \calL_{\sfX} (f)  = \sum_{i=1}^{\dim M} X^i \frac{\de f}{\de x^i}.
\]
where we have given also an explicit expression in local coordinates.
Setting then the Lie bracket of any two functions to vanish
makes $\widetilde\calV^0 \oplus \widetilde\calV^1$ into a GLA.

Then we define the bracket between a vector field $\sfX$ and a
homogeneous element $\sfY_1 \wedge \ldots \wedge \sfY_k \in \widetilde\calV^k$ 
with $k>1$ by the following formula:
\[
\Lie{\sfX}{\sfY_1 \wedge \ldots \wedge \sfY_k} := \sum_{i=1}^{k} (-)^{i +1} 
\Lie{\sfX}{\sfY_i} \wedge \sfY_1\wedge \ldots \wedge \widehat{\sfY}_i  \wedge \ldots \wedge \sfY_k ,
\]
where the bracket on the r.h.s.\ is just the usual bracket on
$\widetilde\calV^1$; we can then extend it to the case of two 
generic multivector fields by requiring it to be linear, graded commutative
and such that for any $\sfX \in  \widetilde\calV^k$, $\ad_{\sfX} := \Lie{\sfX}{\cdot}$
is a derivation of degree $k-1$ w.r.t.\ the wedge product. 

Finally, by iterated application of the Leibniz rule, we can find also an explicit expression  
for the case of a function and a $k$\ndash vector field:
\[
\Lie{\sfX_1 \wedge \cdots \wedge \sfX_k}{f} := \sum_{i=1}^k\;(-)^{k - i} \; \calL_{\sfX_i} (f)\; 
\sfX_1 \wedge \cdots \wedge \widehat\sfX_i \wedge \cdots \wedge \sfX_k
\]
and two homogeneous multivector fields of degree greater than 1:
\[
\begin{aligned}
\Lie{\sfX_1 \wedge \cdots \wedge \sfX_k}{\sfY_1 \wedge \cdots \wedge \sfY_l} & := \\
\sum_{i=1}^{k}\sum_{j=1}^{l}  (-)^{i+j} \;  \Lie{\sfX_i}{\sfY_j}
\wedge \sfX_1\wedge & \cdots \wedge \widehat{\sfX}_i  \wedge  \cdots \wedge \sfX_k 
\wedge \sfY_1\wedge \cdots \wedge \widehat{\sfY}_j  \wedge \cdots \wedge \sfY_l.
\end{aligned}
\]

With the help of these formulae, we can finally check that the 
bracket defined so far satisfies also the Jacobi identity.\footnote{
We give here a sketchy proof; to simplify the notation
the wedge product has not been explicitly written, a small
caret $\miss{\sfV}{i}$ represents the $i$\ndash th component of the missing 
vector field $\sfV$ and $\theta^a_b$ is equal to 1 if $a>b$ and zero otherwise.

Given any three multivector fields $\sfX$, $\sfY$ and $\sfZ$
of positive degree $n$, $l$ and $m$ respectively:
\[
\begin{aligned}
& \Lie{\sfX}{\Lie{\sfY}{\sfZ}} = 
\sum_{i,j}^{l,m} (-)^{i+j}\Lie{\sfX}{\Lie{\sfY_i}{\sfZ_j}
 \; \miss{\sfY}{i} \; \miss{\sfZ}{j} }= \\
& = \sum_{i,j,k}^{l,m,n} (-)^{i+j+k+1} 
[\sfX_k, [\sfY_i, \sfZ_j]]
\; \miss{\sfX}{k}  \; \miss{\sfY}{i} \;  \miss{\sfZ}{j} + 
  \sum_{i,j,k,r \neq i}^{l,m,n} (-)^{i+j + k+r+\theta^r_i}
[\sfX_k, \sfY_r]  [\sfY_i, \sfZ_j]
  \; \miss{\sfX}{k} \; \miss{\sfY}{i,r} \; \miss{\sfZ}{j} + \\
& + \sum_{i,j,k,s \neq j}^{l,m,n} (-)^{i+j+k+s+l-1+\theta^s_l} 
[\sfX_k, \sfZ_s]   [\sfY_i, \sfZ_j]
 \;  \miss{\sfX}{k} \;  \miss{\sfY}{i} \;  \miss{\sfZ}{j,s} = \\
& = \sum_{i,j,k}^{l,m,n} (-)^{i+j+k+1} 
\Big(
[[\sfX_k, \sfY_i], \sfZ_j] 
\; \miss{\sfX}{k}  \; \miss{\sfY}{i} \;  \miss{\sfZ}{j} + 
(-)^{(n+1)(l+1)} [\sfY_i, [\sfX_k, \sfZ_j]]
\; \miss{\sfY}{i} \; \miss{\sfX}{k} \;  \miss{\sfZ}{j} 
\Big)+ \cdots \\
&= \Lie{\Lie{\sfX}{\sfY}}{\sfZ} + (-)^{(n+1)(l+1)} \; \Lie{\sfY}{\Lie{\sfX}{\sfZ}}
\end{aligned}
\]
Analogous computations show that the Jacobi identity is fulfilled
also in the case when one or two of the multivector fields is of degree 0,
while in the case of three functions the identity becomes trivial.
}

This inductive recipe to construct a Lie bracket out of its action on 
the components of lowest degree of the GLA
together with its defining properties completely determines the bracket on 
the whole algebra, as the following proposition summarizes.

\begin{Prop}
There exists a unique extension of the Lie bracket on 
$\widetilde\calV^0 \oplus \widetilde\calV^1$
--- called {\sf Schouten--Nijenhuis bracket} --- 
onto the whole $\widetilde\calV$
\[
\SN{}{}\colon \widetilde\calV^k \otimes  \widetilde\calV^l \to \widetilde\calV^{k+l-1}
\]
for which the following identities hold:

\begin{itemize}
\item[i)] $\SN{\sfX}{\sfY} = - (-)^{(x+1)(y+1)} \; \SN{\sfY}{\sfX}$
\item[ii)] $\SN{\sfX}{\sfY \wedge \sfZ} = \SN{\sfX}{\sfY} \wedge \sfZ 
+ (-)^{(y+1) z} \;\sfY \wedge \SN{\sfX}{\sfZ}$
\item[iii)] $\SN{\sfX}{\SN{\sfY}{\sfZ}} = \SN{\SN{\sfX}{\sfY}}{\sfZ} 
+ (-)^{(x+1)(y+1)} \; \SN{\sfY}{\SN{\sfX}{\sfZ}}$
\end{itemize}

for any triple $\sfX$,$\sfY$ and $\sfZ$ of degree resp. $x$, $y$ and $z$.
\end{Prop}

The sign convention adopted thus far is the original one, as can be 
found for instance in the seminal paper \cite{BFFLS}.
In order to recover the signs we introduced in \ref{def-gla}, we have to shift 
the degree of each element by one, defining the graded Lie algebra of multivector fields
$\calV$ as
\begin{equation} \label{def-calV}
\hspace{2 cm}
\calV := \bigoplus_{i=-1}^{\infty} \calV^i        
\qquad \calV^i := \widetilde\calV^{i+1}\qquad i= -1,0, \ldots,
\end{equation}
which in a shorthand notation is indicated by $\calV := \widetilde\calV [1]$, 
together with the above defined Schouten--Nijenhuis bracket.

The GLA $\calV$ is then turned into a differential graded Lie algebra setting the differential
$\dd \colon \calV \to \calV$ to be identically zero.

\par

We now turn our attention to the particular class of multivector fields
we are most interested in: Poisson bivector fields.
We recall that given a bivector field $\pi \in \calV^1$, we can uniquely define
a bilinear bracket $\Poiss{}{}$ as in \eqref{p-bivector}, which is
by construction skew-symmetric and satisfies Leibniz rule.
The last condition for $\Poiss{}{}$ to be a Poisson bracket
--- the Jacobi identity --- translates into a quadratic equation on the bivector field,
which in local coordinates is:
\[
\begin{array}{c}
\Poiss{\Poiss{f}{g}}{h} + \Poiss{\Poiss{g}{h}}{f} + \Poiss{\Poiss{h}{f}}{g} = 0\\ 
\Updownarrow\\
\pi^{ij} \, \de_j \pi^{kl}  \,\de_j f  \,\de_k g  \,\de_l h +
\pi^{ij}  \,\de_j \pi^{kl}  \,\de_j g  \,\de_k h  \,\de_l f +
\pi^{ij}  \,\de_j \pi^{kl}  \,\de_j h  \,\de_k f  \,\de_l g = 0 \\
\Updownarrow\\
\pi^{ij} \, \de_j \pi^{kl} \; \de_i \wedge \de_k \wedge \de_l =0 
\end{array}
\]
The last line is nothing but the expression in local coordinates 
of the vanishing of the Schouten--Nijenhuis bracket of $\pi$ with itself.
If we finally recall that we defined $\calV$ to be a DGLA with zero differential,
we see that Poisson bivector fields are exactly the solutions to the 
Maurer--Cartan equation \eqref{MC-eq} on $\calV$
\begin{equation} \label{MC-eq-calV}
\dd \pi + \frac{1}{2} \SN{\pi}{\pi} = 0, \qquad \pi \in \calV^1. 
\end{equation}

Finally, formal Poisson structures $\Poiss{}{}_{\ee}$ are associated to a formal 
bivector $\pi \in \ee\, \calV^1\eps$ as in \eqref{formal-poisson} and
the action defined in \eqref{gauge-action} is exactly the gauge group action 
in the sense of Section \ref{dgla}, since the formal diffeomorphisms
acting on $\Poiss{}{}_{\ee}$ are generated by elements of $\calV^0\eps$.

\subsubsection{The DGLA $\calD$} \label{dgla-calD}

The second DGLA that plays a role in the formality theorem is a subalgebra 
of the Hochschild DGLA, whose definition and main properties 
we are going to review in what follows.

To any associative algebra with unit $A$ on a field $\bbK$ 
we can associate the complex of multilinear maps from $A$ to itself.
\[
\calC := \sum_{i=-1}^{\infty} \calC^i \qquad \calC^i := \Hom_\bbK(A^{\otimes(i +1)}, A)
\]
In analogy to what we have done for the case of multivector fields, we shifted 
the degree by one in order to match our convention for the signs that will appear in 
the definition of the bracket.

\begin{figure}[ht] 
\begin{center}
\resizebox{10 cm}{!}{\includegraphics{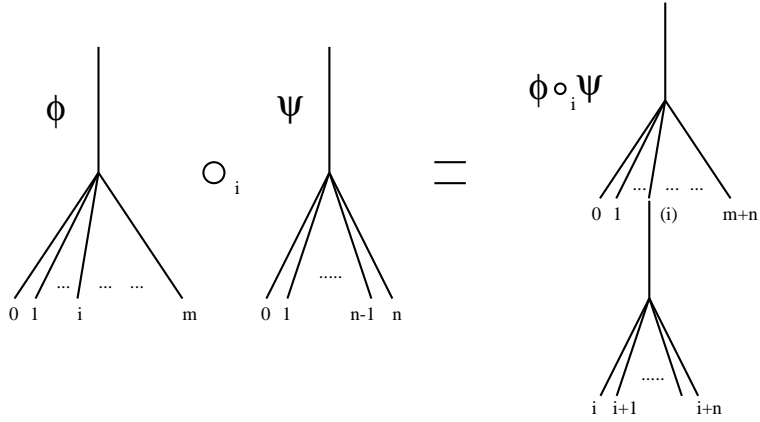}}
\caption{The i- composition.}  \label{hoch1}
\end{center}
\end{figure}

Having the case of linear operators in mind, on which the Lie algebra structure 
arises from the underlying associative structure given by the composition
of operators, we try to extend this notion to multilinear operators.
Clearly, when composing an $(m+1)$\ndash linear operator $\phi$ with 
an $(n+1)$\ndash linear operator $\psi$ we have to specify an inclusion 
$A \hookrightarrow A^{\otimes (m+1)}$ to identify the target space of $\psi$
with one of the component of the domain of $\phi$: loosely 
speaking we have to know where to plug in the output of $\psi$ 
into the inputs of $\phi$. We therefore define a whole family of compositions
$\lbrace \circ_i   \rbrace$ such that for $\phi$ and $\psi$ as above
\[
(\phi \circ_i \psi) (f_0,\dots, f_{m+n}) := \phi( f_0, \ldots,f_{i-1},\psi(f_i,\dots,f_{i+n}), 
f_{i+n+1},\ldots, f_{m+n})
\]
for any $(m+n+1)$\ndash tuple of elements of $A$;
this operation can be better understood through the pictorial representation in Fig. \ref{hoch1}.

We can further sum up with signs all the possible 
partial compositions to find a product 
on $\calC$ --- in fact a pre-Lie structure ---  given by
\[
\phi \circ \psi := \sum_{i=0}^{m} (-)^{n i} \phi \circ_i \psi
\]
with the help of which we can give $\calC$ the structure of a GLA.

\begin{Prop}
The graded vector space $\calC$ together with the Gerstenhaber bracket 
$\Ger{}{}\colon \calC^m \otimes \calC^n \to \calC^{m+n}$ defined 
(on homogeneous elements) by
\begin{equation} \label{def-ger}
\Ger{\phi}{\psi} := \phi \circ \psi - (-)^{m n} \psi \circ \phi
\end{equation}
is a graded Lie algebra, called the Hochschild GLA.
\end{Prop}
\begin{proof}
Since this bracket, introduced by Gerstenhaber in \cite{Ger}, 
is defined as a linear combination
of terms of the form $\phi \circ_i \psi$ and $\psi \circ_i \phi$,
it is clearly linear and homogeneous by construction.
The presence of the sign $(-)^{m n}$ ensures that it is 
also (graded) skew-symmetric, 
since clearly 
\[
\Ger{\phi}{\psi} = 
- (-)^{m n} \Big( \psi \circ \phi -  (-)^{m n} \phi \circ \psi \Big) =
- (-)^{m n} \Ger{\psi}{\phi}
\]
for any $\phi \in \calC^m$ and $\psi \in \calC^n$.

As for the Jacobi identity, we have to prove that
the following holds:
\begin{equation} \label{jacobi-gerst}
\Ger{\phi}{\Ger{\psi}{\chi}} = 
\Ger{\Ger{\phi}{\psi}}{\chi} + (-)^{m n} \Ger{\psi}{\Ger{\phi}{\chi}}
\end{equation}
for any triple $\phi, \psi, \chi$ of multilinear operator 
of degree resp. $m$, $n$ and $p$.
Expanding the first term on r.h.s.\ of \eqref{jacobi-gerst} we get 
\[
\begin{aligned}
&
\Big( \phi \circ \psi - (-)^{m n} \psi \circ \phi \Big) \circ \chi
- (-)^{(m+ n)p}
\chi \circ \Big( \phi \circ \psi - (-)^{m n} \psi \circ \phi \Big) = \\
= &
\sum_{i,k =0}^{m,m+n} \mt{(-)}^{n i + k p} \; (\phi \circ_i \psi ) \circ_k \chi -
\sum_{j,k =0}^{n,m+n} \mt{(-)}^{m(j+n)+ k p} \; (\psi \circ_j \phi) \circ_k \chi +\\
- & 
\sum_{i,k =0}^{m, p} \mt{(-)}^{(m+n)(k+p) + n i} \; \chi \circ_k (\phi \circ_i \psi) +
\sum_{j,k =0}^{n, p} \mt{(-)}^{(m+n)(k+p) + m(j+n)} \; \chi \circ_k (\psi \circ_j \phi)  
\end{aligned}
\]
The first sum can be decomposed according to the following rule
for iterated partial compositions
\[
(\phi \circ_i \psi) \circ_k \chi = 
\left\lbrace\begin{array}{lcc}
(\phi \circ_k  \chi) \circ_i \psi       &\qquad&      k < i\\
\phi \circ_i ( \psi \circ_{k-i} \chi)   &\qquad&    i \leq k \leq i+n\\
(\phi \circ_{k-n}  \chi) \circ_i \psi       &\qquad&      i + n < k
\end{array}\right.
\]
in a term of the form
\[
\mathop{\sum_{i}^{m}}_{ i \leq k \leq i+n}
\mt{(-)}^{n i + kp} \; \phi \circ_i (\psi \circ_{k-i} \chi)=
\sum_{i,k =0}^{m,n} \mt{(-)}^{(n+p) i + k p} \; \phi \circ_i (\psi \circ_k \chi),
\]
whose sign matches the one of the corresponding term coming from 
$(\phi \circ \psi) \circ \chi$ on the l.h.s, plus those 
terms in which the $i$\ndash th and $k$\ndash th composition 
commute, which cancel with the corresponding terms coming from
the expansion of the second term of the r.h.s.\ of \eqref{jacobi-gerst}.

Upon application of the same procedure to the remaining terms, 
the claim follows.
\end{proof}

For a different approach refer to
\cite{St}, where, after having identified
multilinear maps on $A$ with graded coderivations of the free 
cocommutative coalgebra cogenerated by $A$ as a module, the bracket
is interpreted as  the commutator w.r.t.\ the composition of coderivations.

\begin{figure}[hb] 
\begin{center}
\resizebox{8 cm}{!}{\includegraphics{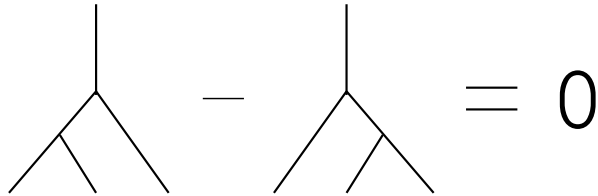}}
\caption{The associativity constraint} \label{assoc}
\end{center}
\end{figure}

Before introducing a differential on $\calC$, we have to 
pick out a particular class of degree one linear operators.
It is clear from the above definitions that associative multiplications
are elements of $\calC^1$ which moreover satisfy the associativity condition.
Writing this equation explicitly in terms of such an element $\Lm$
\begin{equation} \label{ass1}
(f \cdot g )\cdot h =  f \cdot( g \cdot h) 
\Leftrightarrow
\Lm(\Lm(f,g),h) - \Lm(f, \Lm(g,h)) =0 
\end{equation}
we realize immediately that this is --- up to a multiplicative factor ---
the requirement that the Gerstenhaber bracket of $\Lm$ with itself vanishes, since
\begin{equation} \label{ass2}
\begin{aligned}
\Ger{\Lm}{\Lm}(f,g,h) &= \sum_{i=0}^{1} (-)^i  (\Lm \circ_i \Lm )  (f,g,h) - (-)^{1}  
\sum_{i=0}^{1} (-)^i  (\Lm \circ_i \Lm)  (f,g,h) \\
& = 2 \Big( \Lm(\Lm(f,g),h) - \Lm(f, \Lm(g,h)) \Big),
\end{aligned}
\end{equation}
as is shown in a pictorial way in Fig. \ref{assoc}

Now, for each element $\phi$ of degree $k$ of a (DG) Lie algebra $\Lg$,
$\ad_{\phi}:= \Lie{\phi}{}$ is a derivation (of degree $k$),
since the Jacobi identity can also be written as:
\[
\ad_{\phi} \Lie{\psi}{\xi} = \Lie{\ad_{\phi} \psi}{\xi} + (-)^{k m} \Lie{\psi}{\ad_{\phi} \xi} 
\]
for any $\psi \in \Lg^m$ and $\xi \in \Lg^n$.
It is therefore natural to introduce the {\sf Hochschild differential}
\[
\begin{array}{rcl}
\dd_{\Lm} \colon \calC^i &\to& \calC^{i+1}\\
\psi &\mapsto& \dd_{\Lm} \psi := \Ger{\Lm}{\psi}.
\end{array}
\]
The only thing that we still have to check is that $\dd_{\Lm}$ squares to zero, 
which follows immediately from the Jacobi identity and the 
associativity constraint on $\Lm$ expressed in terms of the Gerstenhaber 
bracket as shown in
\eqref{ass1} and \eqref{ass2}:
\[
\begin{aligned}
(\dd_{\Lm} \circ \dd_{\Lm})\, \psi &= \Ger{\Lm}{\Ger{\Lm}{\psi}} = 
\Ger{\Ger{\Lm}{\Lm}}{\psi} - \Ger{\Lm}{\Ger{\Lm}{\psi}} =\\
&= -  \Ger{\Lm}{\Ger{\Lm}{\psi}} 
\qquad \Leftrightarrow \qquad \dd_{\Lm}^2 =0 
\end{aligned}
\]
So we have proved the following
\begin{Prop} \label{def-hoch}
The GLA $\calC$ together with the differential $\dd_{\Lm}$ 
is a differential graded Lie algebra.
\end{Prop}

We can also give an explicit expression of the action of the 
differential on an element $\psi \in \calC^n$:
\[
\begin{aligned}
( \dd_{\Lm} \psi) (f_0, \ldots, f_{n+1}) &=
\sum_{i=0}^{n} (-)^{i+1} \psi (f_0, \ldots, f_{i-1}, f_i \cdot f_{i+1}, \ldots, f_{n+1}) + \\
& + f_0 \cdot \psi (f_1, \ldots, f_{n+1}) +  (-)^{(n+1)} \psi (f_0, \ldots, f_n) \cdot f_{n+1}.
\end{aligned}
\]

As we already mentioned, in the case $A = \Cinf{M}$, 
what we are actually interested in is not the whole 
Hochschild DGLA, but rather a subalgebra of $\calC$: the DGLA of multidifferential 
operators $\widetilde\calD$. It is defined as a (graded) vector space as the collection
$\widetilde\calD := \bigoplus \widetilde\calD^i$ of the subspaces 
$\widetilde\calD^i \subset \calC^i$ consisting of differential operators
acting on smooth functions on $M$.
It is an easy exercise to verify that $\widetilde\calD$ is closed under 
Gerstenhaber bracket and the action of $\dd_{\Lm}$ and thus is a DGL subalgebra.

We stress the fact that $\widetilde\calD$
also includes operators of order $0$,
i.e.\ loosely speaking operators which ``do not differentiate'': this way also 
the associative product $\Lm$ is still an element of $\widetilde\calD^1$.

Having in mind the defining properties of the star product given in Section
\ref{star-prod} and in particular the requirement that 
$B_i(1,f)=0\quad \forall i \in \bbN, f \in \Cinf{M}$,
which ensures that the unity is preserved through deformation, 
we restrict our choice further, considering only differential operators 
which vanish on constant functions; they build a new DGL subalgebra 
$\calD \subset \widetilde\calD$. We remark, however, that $\dd_{\Lm}$
is no longer an inner derivation when restricted to $\calD$, since clearly
the multiplication does not vanish on constants.

Finally, we want to work out also for this DGLA the role played by 
the Maurer--Cartan equation: we will show that in this case this equation
encodes the associativity of the product.

Given an element $\sfB \in \calD^1$, we can interpret $\Lm + \sfB$ as a deformation
of the original product. As shown in \eqref{ass1} and \eqref{ass2}, the associativity
constraint on $\Lm + \sfB$ translates into
\[
\Ger{\Lm+\sfB}{\Lm+\sfB}=0
\]
which in turn, since $\Lm$ is already associative and 
$\Ger{\Lm}{\sfB}=\Ger{\sfB}{\Lm}=\dd_{\Lm} \sfB$
gives exactly the desired Maurer--Cartan equation \eqref{MC-eq}
\begin{equation} \label{MC-eq-calD}
\dd_{\Lm} \sfB + \frac{1}{2} \Ger{\sfB}{\sfB} =0.
\end{equation}

Introducing the formal counterpart of $\calD$, it is clear that 
the deformed product turns out to be nothing but a star product as
in Definition \ref{def-star-prod}, since now $\sfB \in \ee\, \calD^1\eps$
is a formal sum of bidifferential operators. Analogously, the gauge group
is given exactly by formal differential operators and the action on the star 
product is the one given in \eqref{eq-star-prod}, since the adjoint action,
due to the definition of the Gerstenhaber bracket, is nothing but the composition 
of $D_i$ with $B_j$.

\subsection{The first term: $U_1$} \label{hkr-U1}

In this last Section we will give an account for the structures we had 
to introduce and for the two particular cases of DGLA we defined above.

As we already mentioned, our main goal is to prove the formality
of the DGLA $\calD$ of multidifferential operators. 
This approach relies on the existence of a 
previous result by Hochschild, Kostant and Rosenberg \cite{HKR}
which, for any given smooth manifold $M$,  establishes
an isomorphism between the cohomology of the algebra of 
multidifferential operators and the algebra
of multivector fields which, according to our previous definition, 
coincides with its cohomology.
\[
\hkr \colon \calH(\widetilde\calD) \iso \widetilde\calV = \calH(\widetilde\calV )
\]
Actually the original result concerned smooth affine algebraic
varieties, but it can be extended to smooth manifolds, as is shown 
for instance in \cite{Ko2}. 
This isomorphism is induced by the natural map
\[
U_1^{(0)} \colon  \widetilde\calV \longrightarrow \widetilde\calD
\]
which extends the usual identification between vector fields and 
first order differential operators, mapping a homogeneous element 
of the form $\xi_0 \wedge \cdots \wedge \xi_n$ to the multidifferential 
operator whose action on functions $f_0, \ldots , f_n$ 
is given by
\[
\frac{1}{(n+1)!}\; \sum_{\sigma \in S_{n+1}} 
\sgn(\sigma) \; \xi_{\s{0}}(f_0) \cdots \xi_{\s{n}} (f_n),
\]
where we made use of the above mentioned identification
for each $\xi_i$; the definition is extended to $0$\ndash th
order vector fields as the identity map.
Unfortunately this map, which
can be easily checked to be a chain map, fails to preserve the Lie structure,
as can be easily verified already at order $2$. 
Given two homogeneous bivector fields 
$\chi_1 \wedge \chi_2$ and $\xi_1 \wedge \xi_2$,
we can verify explicitly that in general
\[
U_1^{(0)} \left(\Lie{\chi_1 \wedge \chi_2}{\xi_1 \wedge \xi_2} \right)
\neq
\Lie{ U_1^{(0)} (\chi_1 \wedge \chi_2)}{U_1^{(0)}(\xi_1 \wedge \xi_2)}. 
\]
Omitting the subscripts SN and G and the wedge products to ease the notation,
the l.h.s.\ applied to a triple of functions gives
\[
\begin{aligned}
U_1^{(0)} &\left(\Lie{\chi_1}{\xi_1} \, \chi_2 \, \xi_2
- \Lie{\chi_1}{\xi_2} \chi_2 \xi_1
- \Lie{\chi_2}{\xi_1} \chi_1 \xi_2
+ \Lie{\chi_2}{\xi_2} \chi_1 \xi_1
\right) 
(f \tinyotimes g \tinyotimes h) = \\
= \frac{1}{6} &
\Big(\chi_1 \xi_1 f\, \chi_2 g\, \xi_2 h - \xi_1 \chi_1 f\, \chi_2 g\, \xi_2 h
- \chi_1 \xi_2 f\, \chi_2 g\, \xi_1 h + \xi_2 \chi_1 f\, \chi_2 g\, \xi_1 h +\\ 
& -  \chi_2 \xi_1 f\, \chi_1 g\, \xi_2 h + \xi_1 \chi_2 f\, \chi_1 g\, \xi_2 h
+ \chi_2 \xi_2 f\, \chi_1 g\, \xi_1 h + \xi_2 \chi_2 f\, \chi_1 g\, \xi_1 h
\Big) + \text{perm.}
\end{aligned}
\]
while the r.h.s.\ is
\[
\begin{aligned}
& \Lie{\frac{1}{2} \left(\chi_1 \cdot \chi_2 - \chi_2 \cdot \chi_1 \right)}
{\frac{1}{2} \left(\xi_1 \cdot \xi_2 - \xi_2 \cdot \xi_1 \right)}
(f\tinyotimes g \tinyotimes h)=\\
= & \frac{1}{4} 
\Big( \chi_1 (\xi_1 f \, \xi_2 g)\, \chi_2 h + \cdots
\Big).
\end{aligned}
\]
%cr%
However the difference between the two terms is the image
of a closed term in the cohomology of $\calD$. 
We have therefore a way to control the defect 
of this map in being a Lie algebra morphism and 
we can hope to find a way to extend it somehow 
to a morphism whose first order approximation
is this isomorphism of complexes.
This is exactly the role played by the \LL morphism $U$
we will define in the next Sections: in order to give 
a geometric interpretation of this approximation
we will look at the same problem from a dual perspective.

\section{Digression: what happens in the dual} \label{digression}

The whole machinery  of the Kontsevich's construction
can be better understood by looking at the mathematical 
objects and structures we previously introduced from a 
dual point of view.

Given a vector space $V$, polynomials on $V$ can be naturally
identified with symmetric functions on the dual space $V^*$
defining
\[
f(v) := \sum \frac{1}{k!} \; f_k(v \cdots v) \qquad \forall v \in V
 \]
where the coefficients $f_k$ are elements of $\sym{V^*}{k}$.

To extend this construction to the case when $V$ is a 
graded vector space we have to consider the exterior algebra
instead. If we introduce the completion $\Ext{V^*}{}$ of 
this algebra\footnote{To be more precise, we should 
specify the topology w.r.t.\ which we define this completion.
This can be done in a natural way considering $\Sym{V^*}{}$
(resp. $\Ext{V^*}{}$) as the injective limit of the 
$\sym{V^*}{k}$ (resp. $\ext{V^*}{k}$) with the induced 
topology, as in the case of formal power series.},
we can define in a similar way a function in a formal neighborhood 
of $0$ to be given by the formal Taylor expansion in the parameter $\ee$
\[
f(\ee v) := \sum \frac{\ee^k}{k!} \; f_k(v \cdots v) \qquad \forall v \in V.
\]

Following this recipe, a vector field $\sfX$ on $V$ can be identified with
a derivation on $\Ext{V^*}{}$ and Leibniz rule ensures that $\sfX$ is completely
determined by its restriction on $V^*$. In an analogous way 
an algebra homomorphism 
\[
\phi \colon \Ext{W^*}{} \to \Ext{V^*}{},
\]
determines a map $f = \phi^* \colon \Ext{V}{} \to \Ext{W}{}$ 
whose components $f_k$ are completely determined 
by their projection on $W$ as the $\phi_k$ are 
determined by their restriction on $W^*$.

In the following we will need the pointed version of these objects,
namely we will consider the pair $(V, 0)$ as a pointed manifold
and define a (formal) {\sf pointed map} to be an algebra homomorphism
between the reduced symmetric algebras 
(as introduced in \ref{tensor-alg})
\[
\phi \colon \Ext{W^*}{}_{>0} \to \Ext{V^*}{}_{>0},
\]
where the subscript `` $\mt{>0}$'' indicates that we are considering 
the two coalgebras as the (completion of the) quotients 
of $\overline{T}(W^*)$ (resp. $\overline{T}(V^*)$).
Analogously, a {\sf pointed vector field} $\sfX$ is
a vector field which has zero as a fixed point, i.e.\ such that
\[
\sfX(f)(0) = 0 \qquad \forall f
\]
or equivalently such that $(X f)_0 =0$ for every map $f$.

We will further call a pointed vector field 
{\sf cohomological} --- or $Q$\ndash field ---
iff it commutes with itself,
i.e.\ iff $\sfX^2 = \frac{1}{2}[\sfX, \sfX] = 0$ and
{\sf pointed $Q$\ndash manifold} a (formal) pointed 
manifold together with a cohomological vector field.

We turn now our attention to the non commutative case, 
taking a Lie algebra $\Lg$. The bracket 
$\Lie{}{} \colon \Lambda^2 \Lg \to \Lg$ gives rise
to a linear map
\[
\Lie{}{}^* \colon \Lg^* \to \ext{\Lg}{2}^*.
\]
We can extend it to whole exterior algebra to 
\[
\delta \colon \ext{\Lg}{\bullet}^* \to \ext{\Lg}{\bullet+1}^*
\]
requiring that $\delta |_{\Lg^*} \equiv \Lie{}{}^*$ and imposing 
the Leibniz rule to get a derivation.

The exterior algebra can now be interpreted as some odd
analog of a manifold, on which $\delta$ plays the role 
of a (pointed) vector field.
Since the Jacobi identity on $\Lie{}{}$ translates to the equation 
$\delta^2 = 0$, $\delta$ is a cohomological pointed vector field.

If we now consider two Lie algebras $\Lg$ and $\Lh$ 
and endow their exterior algebras with differentials 
$\delta_{\Lg}$ and $\delta_{\Lh}$, 
a Lie algebra homomorphism $\phi \colon \Lg \to \Lh$
will correspond in this case to a chain map $\phi^* \Lh^* \to \Lg^*$,
since 
\[
\phi\left(\Lie{\cdot}{\cdot}_{\Lg}\right) = 
\Lie{\phi(\cdot)}{\phi(\cdot)}_{\Lh}
\qquad \Longleftrightarrow \qquad
\delta_{\Lg} \circ \phi^* = \phi^* \circ \delta_{\Lh}
\]

This is the first glimpse of the correspondence between \LL algebras
and pointed $Q$\ndash manifolds: a Lie algebra is a particular case
of DGLA, which in turn can be endowed with an \LL structure; from
this point of view the map $\phi$ satisfies the same equation
of the first component of an \LL morphism as given in \eqref{FQ-QF} 
for $n=1$.

To get the full picture, we have to extend the previous 
construction to the case of a graded vector space $Z$ which
has odd and even parts. Functions on such a space can be
identified with elements in the tensor product 
$\sym{Z^*}{}:=\sym{V^*}{} \otimes \ext{W^*}{}$, 
where $Z = V \oplus \Pi W$ is the natural decomposition of the 
graded space in even and odd subspaces.\footnote{
In the following we will denote by $\Pi W$ the (odd) space 
defined by a parity reversal on the vector space $W$, which 
can be also written as $W[1]$, using the notation introduced in
Section \ref{dgla}.
}

The conditions for a vector field 
$\delta \colon \sym{Z*}{\bullet} \to \sym{Z*}{\bullet+1}$
to be cohomological can now be expressed in terms of its 
coefficients 
\[
\delta_k \colon  \sym{Z^*}{k} \to \sym{Z^*}{k+1}
\]
expanding the equation $\delta^2 = 0$.
This gives rise to an infinite family of equations:
\[
\left\lbrace
\begin{aligned}
& \delta_0 \, \delta_0 = 0 \\
& \delta_1 \, \delta_0 + \delta_0 \, \delta_1 = 0 \\
& \delta_2 \, \delta_0 + \, \delta_1 \delta_1 + \, \delta_0 \delta_2 =0 \\
& \cdots 
\end{aligned}
\right.
\]

If we now define the dual coefficients 
$m_k := \left(\delta_k|_{Z^*}\right)^*$ and introduce the natural pairing 
$\braket{}{} \colon Z^* \otimes Z \to \bbC$,
we can express the same condition in terms of the maps
\[
m_k \colon \sym{Z}{k+1} \to Z,
\]
paying attention to the signs we have to introduce 
for $\delta$ to be a (graded) derivation.

The first equation ( $m_0 \, m_0 = 0 $ ) tells us that $m_0$ is
a differential on $Z$ and defines therefore a cohomology $\calH_{m_0}(Z)$.

For $k=1$, with an obvious notation, we get
\[
\braket{\delta_1 \, \delta_0 \, f}{x y} = 
\braket{\delta_0 \, f}{m_1( x y)} =
\braket{f}{m_0 ( m_1 ( xy))} 
\]
and
\[
\begin{aligned}
& \braket{\delta_0 \, \delta_1 \, f}{x y} = 
\braket{\delta_1 \, f}{m_0(x)\, y} + (-)^{|x|} \braket{\delta_1 \, f}{x\, m_0(y)} =\\
= & \braket{f}{m_1 ( m_0(x)\, y)} + (-)^{|x|} \braket{f}{m_1 ( x\, m_0(y))}, 
\end{aligned}
\]
i.e.\ $m_0$ is a derivation w.r.t.\ the multiplication defined by $m_1$.

If we now write $Z$ as $\Lg[1]$ and identify the symmetric and
exterior algebras with the d\'ecalage isomorphism 
$\sym{\Lg[1]}{n} \iso \ext{V[n]}{n}$,
$m_1$ can be interpreted as a bilinear skew-symmetric operator on $\Lg$.

The next equation, which involves $m_1$ composed with itself, tells us
exactly that this operator is indeed a Lie bracket for which the 
Jacobi identity is satisfied up to terms containing $m_0$, 
i.e.\ --- since $m_0$ is a differential --- up to homotopy. 
 
Putting the equations together, this gives rise to a strong homotopy 
Lie algebra structure on $\Lg$, thus establishing a one-to-one
correspondence between pointed $Q$\ndash manifolds and SHLA's, 
which in turn are equivalent to \LL algebras, 
as we already observed in Section \ref{dgla}.\\
 
Finally, to complete this equivalence and to express 
the formality condition \eqref{FQ-QF} more explicitly,
we spell out the equations for the coefficients of a 
{\sf $Q$\ndash map}, i.e.\ a (formal) pointed map between
two $Q$\ndash manifolds $Z$ and $\tZ$ which commutes 
with the $Q$\ndash fields; namely:
\begin{equation} \label{FQ-QF2}
\begin{array}{rcl}
\phi \colon \sym{\tZ^*_{>0}}{} &\longrightarrow& \sym{Z^*_{>0}}{} \\
\, &\text{s. t.}& \, \\
\phi \circ \tilde{\delta} &=& \delta \circ \phi.
\end{array}
\end{equation}

As for the case of the vector field $\delta$, we consider only the restriction
of this map to the original space $\tZ$ and define the coefficients
of the dual map as
\[
U_k := \left( \phi_k|_{\tZ^*}\right)^* \colon \sym{Z}{k} \to \tZ.
\]

With the same notation as above, we can express the condition
\eqref{FQ-QF2} on the dual coefficients with the help of the 
natural pairing. 
The first equation reads:
\[
\begin{array}{rcl}
\braket{\phi \, \tilde{\delta} \, f}{x} &=& 
\braket{\delta \, \phi \, f}{x}\\
\, &\Downarrow& \,\\
\braket{\tilde{\delta_0} \, f}{ U_1(x)} &=&
\braket{\phi \, f}{m_0(x)}\\ 
\, &\Downarrow& \,\\
\braket{f}{\widetilde{m}_0(U_1(x))} &=& 
\braket{f}{U_1(m_0(x))}.
\end{array}
\]

As we could have guessed from the discussion in Section \ref{dgla},
the first coefficient $U_1$ is a chain map w.r.t.\ the differential
defined by the first coefficient of the $Q$\ndash structures.
\[
\left[ U_1 \right] \colon \calH_{m_0}(Z) \to \calH_{\widetilde{m}_0}(\tZ).
\]
An analogous computation gives the equation for
the next coefficient:
\[
\widetilde{m}_1 (U_1(x)\, U_1(y)) + \widetilde{m}_1 (U_2 (x\,y)) =
U_2( m_0(x)\,y) + (-)^{|x|} U_2(x\, m_0(y)) + U_1(m_1(x\,y)),
\]
which shows that $U_1$ preserves the Lie structure induced 
by $m_1$ and $\widetilde{m}_1$ up to terms containing $m_0$ and 
$\widetilde{m}_0$, i.e.\ up to homotopy.

This is exactly what we were looking for: as the map 
$U_1^{(0)}$ defined in Section \ref{hkr-U1} is a chain map
which fails to be a DGLA morphism, a $Q$\ndash map $U$
(or equivalently an \LL morphism) induces a map $U_1$
which shares the same property. 

We restrict thus our attention to DGLA's, considering
now a pair of pointed $Q$\ndash manifolds $Z$ and $\tZ$
such that $m_k = \widetilde{m}_k = 0$ for $k>1$. 
Equivalently, we consider two \LL algebras as in 
Example \ref{DGLA-LL}, whose coderivation have
only two non-vanishing components.

A straightforward computation which follows the same steps 
as above for $k=1,2$, leads in this case to the following
condition on the $n$\ndash th coefficient of $U$:
\begin{equation} \label{FQ-QF3}
\begin{aligned}
\widetilde{m}_0  \left(U_n (x_1 \cdots x_n)\right) 
+ \frac{1}{2}  \sum_{\bisum{I \sqcup J = \lbrace 1, \ldots n \rbrace}{I,J \neq \emptyset}} 
\!\!\!\!\! &  \eee_{x} (I, J) \; \widetilde{m}_1\left( U_{|I|}(x_{I}) \cdot U_{|J|}(x_{J}) \right) =\\
= \sum_{k = 1}^{n} \;\; & \eee_{x}^k \; 
 U_{n}\left( m_0 (x_k) \cdot x_1 \cdots \widehat{x}_k \cdots x_n \right) +\\
+ \frac{1}{2} \sum_{k \neq l} \;\; & \eee_{x}^{k l} \; 
U_{n-1} \left( m_1(x_k \cdot x_l)\cdot x_1 \cdots 
\widehat{x}_k \cdots \widehat{x}_l  \cdots x_n\right) 
\end{aligned}
\end{equation}

To avoid a cumbersome expression involving lots of signs,
we introduced a shorthand notation $\eee_x (I, J)$ for the 
Koszul sign associated to the $(|I|, |J|)$\ndash shuffle permutation 
associated to the partition $I \sqcup J = \{1, \ldots , n\}$\footnote{
Whenever a vector space $V$ is endowed with
a graded commutative product, the {\sf Koszul sign} $\eee(\sigma)$
of a permutation $\sigma$ is the sign defined by
\[
x_1 \cdots x_n = \eee(\sigma) \; x_{\s{1}} \cdots x_{\s{n}} 
\qquad x_i \in V.
\]
An $(l, n-l)$\ndash shuffle permutation is a permutation $\sigma$
of $(1, \ldots, n)$ such that $\sigma(1) < \cdots < \sigma(l)$
and  $\sigma(l+1) < \cdots \sigma(n)$. The shuffle permutation
associated to a partition 
$I_1 \sqcup \cdots \sqcup I_k = \{1 , \ldots , n\}$ 
is the permutation that takes first all the elements indexed 
by the subset $I_1$ in the given order, then those indexed
by $I_2$ and so on.
}and $\eee_x^k$  (resp. $\eee_x^{k l}$)
for the particular case $I = \{k\}$ (resp. $I =\{ k, l\}$);   
we further simplified the expression adopting 
the multiindex notation $x_{I} := \prod_{i \in I} x_i$.

This expression will be specialized in next Section to the case 
of the \LL morphism introduced by Kontsevich to give a formula for 
the star product on $\bbR^d$: we will choose as $Z$ the DGLA $\calV$
of multivector fields and as $\tZ$ the DGLA $\calV$ of
multidifferential operators and derive the equation that the
coefficients $U_n$ must satisfy to determine the required 
formality map.\\

As a concluding act of this digression, we will establish
once and for all the relation between the formality of $\calD$
and the solution of the problem of classifying all possible
star products on $\bbR^d$.

As we already worked out in Section \ref{dgla},
the associativity of the star product as well as 
the Jacobi identity for a bivector field are encoded
in the Maurer--Cartan equations \ref{MC-eq-calD} resp.
\ref{MC-eq-calV}. In order to translate these equations
in the language of pointed $Q$\ndash manifolds,
we have first to introduce the generalized Maurer--Cartan 
equation on an (formal) \LL algebra $(\Lg\eps, Q)$:
\[
Q( \exp \ee\, x) = 0 \qquad x \in \Lg^1\eps,
\]
where the exponential function $\exp$ maps an element of degree
$1$ to a formal power series in $\ee \Lg\eps$.

From a dual point of view, this amounts to the request that 
$x$ is a fixed point of the cohomological vector field $\delta$,
i.e.\ that for every $f$ in $\sym{\Lg^*\eps[1]}{}$
\[
\delta \, f( \ee \, x) =0. 
\]

Since $(\delta f)_k = \delta_{k-1} f$, expanding the previous
equation in a formal Taylor series and using the pairing as 
above to get $\braket{\delta_{k-1} f}{x \cdots x} = 
\braket{f}{m_{k-1}(x \cdots x)}$, the generalized Maurer--Cartan
equation can be written in the form
\begin{equation} \label{gen-MC-eq}
\sum_{k=1}^{\infty} \frac{\ee^k}{k!} m_{k-1}(x \cdots x) = 
\ee \; m_0 (x) + \frac{\ee^2}{2}\; m_1(x\,x) + o(\ee^3) = 0.
\end{equation}

It is evident that (the formal counterpart of) equation \ref{MC-eq}
is recovered as a particular case when $m_k = 0$ for $k>1$.

Finally, as a morphism of DGLA's preserves the solutions 
of the Maurer--Cartan equation, since it commutes both
with the differential and with the Lie bracket,
an \LL morphism 
$\phi \colon \sym{(\Lh^*\eps[1])}{} \to \sym{(\Lg^*\eps[1])}{}$, 
according to \eqref{FQ-QF2},
preserves the solutions of the above generalization;
with the usual notation, if $x$ is a solution to 
\eqref{gen-MC-eq} on $\Lg\eps$, 
\[
U(\ee \, x) = \sum_{k=1} \frac{\ee^k}{k!} \; U_k (x \cdots x)
\]
is a solution of the same equation on $\Lh$.

The action of the gauge group on the set $\maurer(\Lg)$
can analogously be generalized to the case of \LL algebras
and a similar computation shows that, if $x$ and $x'$ 
are equivalent modulo this generalized action, 
their images under $U$ are still equivalent solutions.\\

In conclusion, reducing the previous discussion 
to the specific case we are interested in, namely when
$\Lg=\calV$ and $\Lh=\calD$, given an \LL morphism $U$
we have a formula to construct out of any (formal)
Poisson bivector field $\pi$ an associative star product 
given by
\begin{equation} \label{U-pi}
U(\pi) = \sum_{k=0} \frac{\ee^k}{k!} \; U_k (\pi \cdots \pi)
\end{equation}
where we reinserted the coefficient of order $0$ 
corresponding to the original non deformed product.
If moreover $U$ is a quasi-isomorphism, the correspondence
between (formal) Poisson structures on $M$ and formal deformations
of the pointwise product on $\Cinf{M}$ is one-to-one: 
in other terms once we give a formality map, we have solved 
the problem of existence and classification of star products
on $M$.\\

This is exactly the procedure followed by Kontsevich 
to give his formula for the star product on $\bbR^d$.

\section{The Kontsevich formula}

In this Section we will finally give an explicit expression of 
Kontsevich's formality map from $\calV$ to $\calD$ which induces 
the one-to-one map from (formal) Poisson structures on $\bbR^d$
to star products on $\Cinf{\bbR^d}$. 

The main idea is to introduce a pictorial way to describe how 
a multivector field can be interpreted as a multidifferential operator
and to rewrite the equations introduced in \ref{FQ-QF2}
in terms of graphs.

As a toy model we can consider the Moyal star product introduced in
Section \ref{star-prod} and give a pictorial version of formula 
\eqref{gen-moyal} as follows:

\begin{figure}[ht] 
\begin{center}
\resizebox{12 cm}{!}{\includegraphics{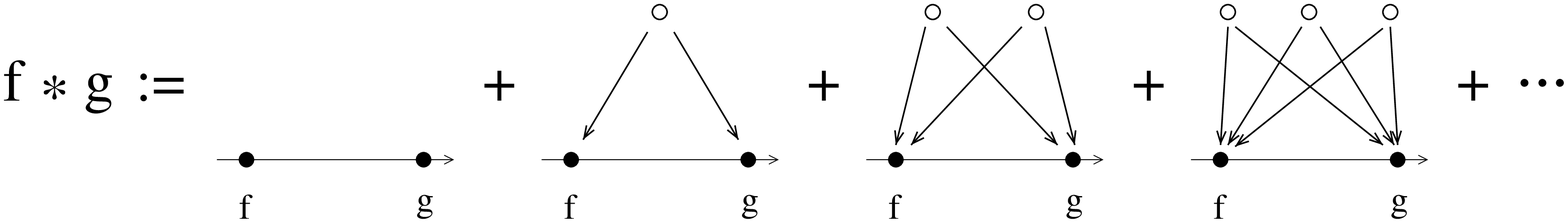}}
\caption{A pictorial representation of the first terms of the Moyal star
  product.} \label{moyal}
\end{center}
\end{figure}

To the $n$\ndash th term of the series we associate a graph with $n$
``unfilled'' vertices -- which represent the $n$ copies of the Poisson 
tensor $\pi$ --  and two ``filled'' vertices -- which stand for the two 
functions that are to be differentiated; the left (resp. right) arrow emerging
from the vertex corresponding to $\pi^{ij}$ represent  $\de_i$ (resp.  
$\de_j$)  acting on $f$ (resp. $g$) and the sum over all indices involved
is understood.

This setting can be generalized introducing vertices of higher order, 
i.e.\ with more outgoing arrows, to represent multivector fields and 
letting arrows point also to ``unfilled'' vertices, to represent the
composition of differential operators: in the Moyal case, since the 
Poisson tensor is constant such graphs do not appear.

The main intuition behind the Kontsevich formula for the star product
is that one can introduce an appropriate set of graphs and assign to
each graph $\Gamma$ a multidifferential operator $B_{\Gamma}$
and a weight $w_{\Gamma}$ in such a way that the map that sends
an $n$\ndash tuple of multivector fields to the corresponding 
weighted sum over all possible graphs in this set 
of multidifferential operators is an \LL morphism.

This procedure will become more explicit in the next Section, 
where we will go into the details of Kontsevich's construction.

\subsection{Admissible graphs, weights and $B_\Gamma$'s}

First of all, we have to introduce the above mentioned 
set of graphs we will deal with in the following. 

\begin{Def} \label{def-adm-graphs}
The set $\calG_{n,\bar{n}}$ of {\sf admissible} graphs
consists of all connected graphs $\Gamma$ which satisfy the following properties:
\begin{itemize}
\item[-] the set of vertices $V(\Gamma)$ is decomposed in two ordered subsets 
$V_1(\Gamma)$ and $V_2(\Gamma)$ 
isomorphic to $\lbrace 1, \ldots , n\rbrace$ resp. 
$\lbrace \bar{1} , \ldots , \bar{n}\rbrace$ whose elements are 
called vertices of the first resp. second type;
\item[-] the following inequalities involving the number of vertices of the two 
types are fulfilled: 
$n \geq 0$\/, $\bar{n} \geq 0$ and $2 n + \bar{n} - 2 \geq 0$;
\item[-] the set of edges $E(\Gamma)$ is finite and does not contain
{\sf small loops}, i.e.\ edges starting and ending at the same vertex;
\item[-] all edges in $E(\Gamma)$ 
are oriented and start from a vertex of the first type;
\item[-] the set of edges starting at a given vertex $v \in V_1(\Gamma)$,
which will be denoted in the following by $\Star(v)$, is ordered. 
\end{itemize}
\end{Def}

\begin{Exa} {\bf Admissible graphs}

Graphs $i)$ and $ii)$ in Fig. \ref{adm-graphs} are admissible, while
graphs $iii)$ and $iv)$ are not.

\begin{figure}[ht] 
\begin{center}
\resizebox{12 cm}{!}{\includegraphics{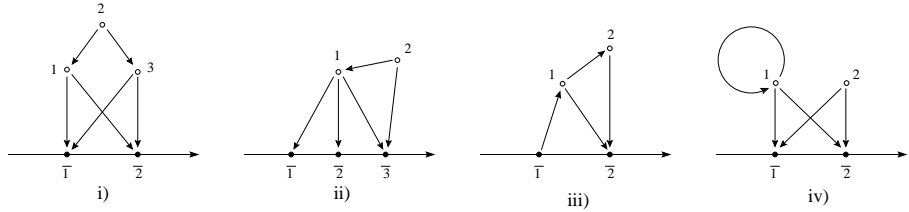}}
\caption{Some examples of admissible and non-admissible graphs.}  \label{adm-graphs}
\end{center}
\end{figure}

\end{Exa}

We now introduce the procedure to associate to each pair 
$(\Gamma, \xi_1 \otimes \cdots \otimes \xi_n)$ consisting of a graph 
$\Gamma \in \calG_{n,\bar{n}}$ with $2 n + m - 2$ edges and of a tensor 
product of $n$ multivector fields on $\bbR^d$
a multidifferential operator $B_{\Gamma} \in \calD^{\bar{n}-1}$.

\begin{itemize}
\item We associate to each vertex $v$ of the first type with $k$
outgoing arrows the skew-symmetric tensor $\xi_i^{j_1, \ldots, j_k}$
corresponding to a given $\xi_i$ via the natural identification.
\item We place a function at each vertex of the second type.
\item We associate to the $l$\ndash th arrow in $\Star(v)$ 
a partial derivative w.r.t.\ the coordinate 
labeled by the $l$\ndash th index of $\xi_i$ acting on
the function or the tensor appearing at its endpoint. 
\item We multiply such elements in the order prescribed 
by the labeling of the graph.
\end{itemize}

As an example, the multidifferential operator corresponding
to the first  graph in Fig.\ref{adm-graphs} and to the triple 
$(\alpha, \beta, \gamma)$ of bivector fields is given by
\[
U_{\Gamma_1}(\alpha, \beta,\gamma)(f,g) :=
\beta^{b_1 b_2}\, \de_{b_1} \alpha^{a_1 a_2} \,
\de_{b_2} \gamma^{c_1 c_2}\, \de_{a_1} \de_{c_1} f \, \de_{a_2} \de_{c_2}
g,
\]
while the operator corresponding to the second graph and 
the pair $(\pi, \rho)$is 
\[
U_{\Gamma_1}(\pi, \rho)(f,g,h) :=
\pi^{p_1 p_2} \de_{p_1} \rho^{r_1 r_2 r_3} 
\de_{r_1} f \de_{r_2} g \de_{r_3} \de_{p_2} h
\]

This construction gives rise for each $\Gamma$ to a linear map 
$U_{\Gamma} \colon T^n(\calV) \to \calD$ which is equivariant w.r.t.\
the action of the symmetric group, i.e.\ permuting the order in which 
we choose the edges we get a sign 
equal to the signature of the permutation.
The main point in Kontsevich's formality theorem was to show that there 
exist a choice of weights $w_{\Gamma}$ such that the linear combination
\[
U := \sum_{\Gamma} w_{\Gamma} B_{\Gamma}
\]
defines an \LL morphism, where the sum runs over all admissible graphs.

These weights are given by the product of a combinatorial coefficient 
times the integral of a differential form $\omega_{\Gamma}$ 
over the configuration space $C_{n,\bar{n}}$ defined in the following.
The expression of the weight $w_{\Gamma}$ associated to 
$\Gamma \in \calG_{n,\bar{n}}$ is then:
\begin{equation} \label{def-weight}
w_{\Gamma} := \prod_{k=1}^{n} \frac{1}{(\# \Star(k))!}
\frac{1}{(2 \pi)^{2 n + \bar{n} -2}} 
\int_{\bar{C}^{+}_{n,\bar{n}}} \!\!\!\!\!\omega_{\Gamma}
\end{equation}
if $\Gamma$ has exactly $2 n + \bar{n} - 2 $ edges, while the weight
is set to vanish otherwise.
The definition of $\omega_{\Gamma}$ and of the configuration space can be better
understood if we imagine embedding the graph $\Gamma$ in the upper half plane
$\calH := \{z \in \bbC |\; \Im(z) \geq 0\}$ binding the vertices of the second type 
to the real line.

We can now introduce the open configuration space of the $n+\bar{n}$ 
distinct vertices of $\Gamma$ as the smooth manifold:
\[
\begin{aligned}
\Conf_{n,\bar{n}}:=\Big\lbrace (z_1, \ldots , z_n, z_{\bar{1}}, \ldots, z_{\bar{n}})
\in \bbC^{n+\bar{n}} \Big|\;&
z_i \in \calH^{+}, z_{\bar{i}} \in \bbR, \\
\, &  z_i\neq z_j \;\text{for}\; i\neq j,\; 
z_{\bar{i}} \neq z_{\bar{j}} \;\text{for}\; \bar{i} \neq \bar{j} \Big\rbrace.
\end{aligned}
\]
In order to get the right configuration space we have to quotient 
$\Conf_{n,\bar{n}}$ by the action of the $2$\ndash dimensional
Lie group $G$ consisting of translations in the horizontal direction
and rescaling, whose action on a given point $z \in \calH$ is given by:
\[
\qquad\qquad z \mapsto a z + b \qquad \qquad a \in \bbR^{+}, b\in \bbR.
\]
In virtue of the condition imposed on the number of vertices in
\eqref{def-adm-graphs}, the action of $G$ is free; therefore the quotient
space, which will be denoted by $C_{n,\bar{n}}$, is again a smooth manifold,
of (real) dimension $2 n + \bar{n} -2$.

Particular care has to be devoted to the case when the graph has no 
vertices of the second type. In this situation, having no points on the real line,
the open configuration space can be defined as a subset of $\bbC^{n}$ 
instead of $\calH^{n}$ and we can introduce a more general Lie group $G'$,
acting by rescaling and translation in any direction; the quotient space 
$C_{n}:= \Conf_{n,0}/_{G'}$ for $n \geq 2$ is again a smooth manifold,
of dimension $2n-3$.

In order to get a connected manifold, we restrict further our attention to 
the component $C^{+}_{n,\bar{n}}$ in which the vertices of the second type
are ordered along the real line in ascending order, namely:
\[
C^{+}_{n,\bar{n}} := \Big\lbrace  (z_1, \ldots , z_n, z_{\bar{1}}, \ldots, z_{\bar{n}})
\in C_{n,\bar{n}} \Big|\;
z_{\bar{i}} < z_{\bar{j}} \;\text{for}\; \bar{i} < \bar{j} \Big\rbrace.
\]

\begin{figure}[ht] 
\begin{center}
\resizebox{6 cm}{!}{\includegraphics{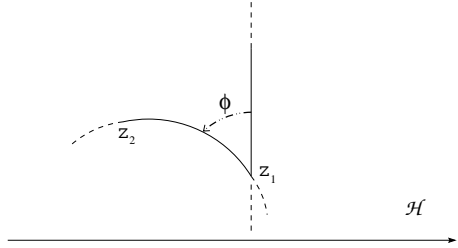}}
\caption{The angle map $\phi$}  \label{angle}
\end{center}
\end{figure}

On these spaces we can finally introduce the differential form
$\omega_{\Gamma}$. We first define an {\sf angle map}
\[
\phi \colon C_{2,0} \longrightarrow S^{1}
\] 
which associates to each pair of distinct points $z_1, z_2$ in the upper half
plane the angle between the geodesics w.r.t.\ the Poincar\'e metric
connecting $z_1$ to $+\,i\, \infty$ and to $z_2$, measured in the
counterclockwise direction (cfr. Fig. \ref{angle}).

The differential of this function is now a well-defined $1$\ndash form
on $C_{2,0}$ which we can pull-back to the configuration space
corresponding to the whole graph with the help of the natural
projection $\pi_{e}$ associated to each edge $e=(z_i, z_j)$ of $\Gamma$
\[
\begin{array}{rccc}
\pi_{e} \colon & C_{n,\bar{n}} & \longrightarrow & C_{2, 0}\\
\, & (z_1, \ldots , z_{\bar{n}})   & \mapsto         & (z_i, z_j) 
\end{array}
\] 
to obtain $d \phi_{e} := \pi_{e}^{*}\, d \phi 
\in \Omega^1 (C_{n,\bar{n}} )$.
The form that appears in the definition of the weight $w_{\Gamma}$ can
now be defined as
\[
\omega_{\Gamma} := \bigwedge_{e \in \Gamma} d \phi_{e}
\]  
where the ordering of the $1$\ndash forms in the product is the one
induced on the set of all edges by the ordering on the (first)
vertices and the ordering on the set $\Star(v)$ of edges emerging from
the vertex $v$. We want to remark hereby that, as long as we consider
graphs with $2 n + \bar{n} -2$ edges, the degree of the form matches exactly
the dimension of the space over which it has to be integrated,
which gives us a real valued weight.

This geometric construction has a more natural interpretation if one 
derives the Kontsevich formula for the star product from a path integral
approach, as it was done for the first time in \cite{CF1}.

For the weights to be well-defined, we also have to require that the
integrals involved converge. However, as the geometric construction 
of $\phi$ suggests, as soon as two points approach each
other, the differential form $d \phi$ is not defined. The solution to this
problem has already been given implicitly in \eqref{def-weight}: the
differential form is not integrated over the open configuration space,
but on a suitable compact space whose definition and properties are
contained in the following
\begin{Lem}
For any configuration space $C_{n,\bar{n}} $ (resp. $C_{n}$) there 
exists a compact space $\bar{C}_{n,\bar{n}} $ (resp. $\bar{C}_{n}$) 
whose interior is the open configuration space and 
such that the projections $\pi_{e}$, the angle map $\phi$ and thus the 
differential form $\omega_{\Gamma}$ extend smoothly to the 
corresponding compactifications.
\end{Lem}
The compactified configuration spaces are (compact) smooth manifolds
with corners. We recall that a smooth manifold with corner of
dimension $m$ is a topological Hausdorff space $M$ which is locally
homeomorphic to $\bbR^{m-n} \times \bbR_{+}^n$ with $n=0,\ldots,m$. 
The points $x \in M$ whose local
expression in some (and thus any) chart has the form 
$x_1, \ldots , x_{m-n}, 0, \ldots ,0)$ 
are said to be {\sf of type $n$} and form submanifolds of $M$ called
{\sf strata} of codimension $n$.

The general idea behind such a compactification is that the naive approach
of considering the closure of the open space in the cartesian product
would not take into account the different speeds with which two or 
more points ``collapse'' together on the boundary of the configuration 
space.

For a more detailed description of the compactification we refer the
reader to \cite{FMP} for an algebraic approach
and to \cite{AS} and \cite{BT} for an explicit description
in local coordinates. More recently Sinha \cite{S} gave a
simplified construction in the spirit of Kontsevich's original ideas.
In \cite{AMM}
the orientation of such spaces and of their codimension one strata --
whose relevance will be clarified in the following -- is discussed.

Finally, the integral in \eqref{def-weight} is well-defined and yields
a weight $w_{\Gamma} \in \bbR$ for any admissible graph $\Gamma$,
since we defined $w_{\Gamma}$ to be non zero only when $\Gamma$ 
has exactly $2 n + m -2$ edges, i.e.\ when the degree of
$\omega_{\Gamma}$ matches the dimension of the 
corresponding configuration space.

\subsection{The proof: Lemmas, Stokes' theorem, Vanishing theorems}

Having defined all the tools we will need, we can now give a sketch
of the proof.

In order to verify that $U$ defines the required \LL morphism
we have to check that the following conditions hold:
\begin{itemize}
\item[I] The first component of the restriction of $U$ to $\calV$ is -- up
to a shift in the degrees of the two DGLAs -- the natural map introduced in 
Section \eqref{hkr-U1}.
\item[II] $U$ is a graded linear map of degree $0$.
\item[III] $U$ satisfies the equations for an \LL morphism 
defined in Section \eqref{digression}.
\end{itemize}

\begin{Lem} \label{lem-I} {\bf I}
The map
\[
U_1 \colon \calV \longrightarrow \calD
\]
is the natural map that identifies each multivector field with the
corresponding multiderivation.
\end{Lem}
\begin{proof}
The set $\calG_{1,\bar{n}}$ consist of only one element, namely the graph 
$\Gamma_{\bar{n}}$
with one vertex of the first type with $2\cdot 1 + \bar{n} - 2 =\bar{n}$ 
arrows with an equal number of vertices of the second type as endpoints.

\begin{figure}[ht!] 
\begin{center}
\resizebox{8 cm}{!}{\includegraphics{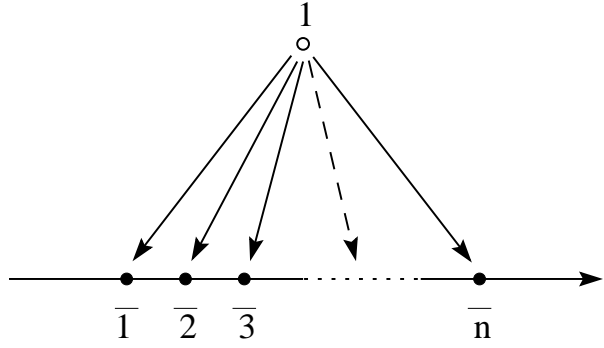}}
\caption{The admissible graph $\Gamma_{\bar{n}}$}  \label{gamma-bar-n}
\end{center}
\end{figure}
To each $k$\ndash vector field $\xi$ we associate thus the multidifferential 
operator given by
\[
U_{\Gamma_{\bar{n}}}(\xi)(f_{\bar{1}}, \ldots, f_{\bar{n}}):=
w_{\Gamma_{\bar{n}}}\; \xi^{i_{\bar{1}}, \ldots, i_{\bar{n}}}\; 
\de_{i_{\bar{1}}} f_{\bar{1}} \cdots \de_{i_{\bar{n}}} f_{\bar{n}}.
\]
An easy computation shows that the integral of $\omega_{\Gamma_{\bar{n}}}$
over $\bar{C}_{1,\bar{n}}$ cancels the power of $\frac{1}{2 \pi}$ and
leaves us with the right weight 
\[
w_{\Gamma_{\bar{n}}} = \frac{1}{\bar{n}!} 
\]
we expect for $U_1$ to be the natural map that induces the HKR isomorphism.
\end{proof}

\begin{Lem} \label{lem-II} {\bf II}
The $n$\ndash th component 
\[
U_n := \sum_{\bar{n} =1}^{\infty} 
\sum_{\Gamma \in \calG_{n,\bar{n}}} w_{\Gamma} B_{\Gamma}
\]
has the right degree for $U$ to be an \LL morphism.
\end{Lem}
\begin{proof}
To each vertex $v_i$ with $\# \Star(v_i)$ outgoing arrows
corresponds an element of $\calV^{r_i} = \widetilde\calV^{r_i +1}$
where $r_i = \# \Star(v_i)$. On the other side, each graph
with $\bar{n}$ vertices of the second type together with
an $n$\ndash tuple of multivector fields gives rise to a differential 
operator of degree $s = \bar{n} -1$.
Since we consider only graphs with $2 n + \bar{n} -2$ edges and this
is equal by construction to 
\[
\sum_{i=1}^{n} \# \Star(v_i),
\]
the degree of $U_n (\xi_1, \ldots, \xi_n)$ can be written as
\[
s = (2 n + \bar{n} -2 ) +1 -n = \sum_{i=1}^n r_i +1 - n
\]
which is exactly the prescribed degree 
for the $n$\ndash th component of an \LL morphism.

Although the construction we gave in the previous section
involves a tensor product of multivector fields, 
the signs and weights in $U_n$ are chosen in such a
way that, upon symmetrization, it descends to the symmetric algebra.
\end{proof}

We come now to the main part of Kontsevich's construction: 
the geometric proof of the formality.

First of all we have to extend our morphism $U$ to 
include also a $0$\ndash th component which 
represents the usual multiplication between smooth functions
--- the associative product we want to deform via the higher
order corrections. 
We can now specialize the $L_{\infty}$ condition \eqref{FQ-QF3} 
to the case at hand, where $m_0$ $\widetilde{m}_0$ 
can be expressed in 
terms of of the Taylor coefficients $U_n$ as:

\begin{equation} \label{form}
\begin{aligned}
\sum_{l=0}^{n} & \sum_{k=-1}^{m} \sum_{i=0}^{m-k} \eee_{kim}\!\!\!\!
\sum_{\sigma \in S_{l, n-l}} \!\!\!\!\!\eee_{\xi}(\sigma)\;
U_l \Big(\xi_{\s{1}}, \ldots , \xi_{\s{l}}\Big)\\
&\Big(f_0 \tinyotimes \cdots \tinyotimes f_{i-1} \tinyotimes 
U_{n-l} ( \xi_{\s{l+1}}, \ldots, \xi_{\s{n}} ) 
(f_i \tinyotimes \cdots \tinyotimes f_{i+k}) \tinyotimes f_{i+k+1}  
\tinyotimes \cdots \tinyotimes f_m\Big)\\  
=& \sum_{i\neq j =1}^{n} \eee_{\xi}^{ij}\; U_{n-1} (\xi_i \circ \xi_j, \xi_1,
\ldots, \widehat{\xi_i},\ldots, \widehat{\xi_j},\ldots, \xi_n )
(f_0 \tinyotimes \cdots \tinyotimes f_n),
\end{aligned}
\end{equation}

where
\begin{itemize}
\item[-] $\{\xi_j\}_{j=1,\ldots, n}$ are multivector fields;
\item[-] $f_0, \ldots, f_m$ are the smooth functions on which the multidifferential 
operator is acting;
\item[-] $S_{l, n-l}$ is the subset of $S_n$ consisting 
of $(l, n-l)$\ndash shuffles
\item[-] the product $\xi_i \circ \xi_j$ is defined in such a way that
the Schouten--Nijenhuis bracket can be expressed in terms of this 
composition by a formula similar to the one relating
the Gerstenhaber bracket to the analogous composition $\circ$
on $\calD$ given in \ref{dgla-calD};
\item[-] the signs involved are defined as follows: $\eee_{kim} := (-1)^{k (m+i)}$,
$\eee_{\xi}(\sigma)$ 
is the Koszul sign associated to the permutation $\sigma$
and $\eee_\xi^{ij}$ is defined as in \eqref{FQ-QF3}.
\end{itemize}
This equation encodes the formality condition since the l.h.s.\
corresponds to the Gerstenhaber bracket between multidifferential 
operators while the r.h.s.\ contains ``one half'' of the
Schouten--Nijenhuis bracket; the differentials do not appear 
explicitly since on $\calV$ we defined $\dd$ to be 
identically zero, while on $\calD$ it is expressed in terms of 
the bracket with the multiplication $\Lm$, which we included in 
the equation as $U_0$.

For a detailed explanation of the signs involved we refer once more to 
\cite{AMM}.\\

We can now rewrite equation \eqref{form} in a form 
that involves again admissible graphs and weights to show 
that it actually holds. It should be clear from the 
previous construction of the coefficients $U_k$ that
the difference between the l.h.s.\ and the r.h.s.\ of 
equation \eqref{form} can be written as a linear combination
of the form 
\begin{equation} \label{form2}
\sum_{\Gamma \in \calG_{n,\bar{n}}} c_{\Gamma} 
U_{\Gamma}(\xi_1, \ldots, \xi_n)(f_0 \tinyotimes \cdots \tinyotimes f_n)
\end{equation}
where the the sum runs in this case over the set of admissible 
graphs with $2 n + \bar{n} -3$ edges. Equation \eqref{form}
is thus fulfilled for every $n$ if these coefficients $c_{\Gamma}$
vanish for every such graph.

The main tool to prove the vanishing of these coefficients 
is the Stokes Theorem for manifolds with corners, which 
ensures that also in this case the integral of an exact form 
$\dd \Omega$ on a manifold $M$ can be expressed as the integral
of $\Omega$ on the boundary $\de M$. 
In the case at hand, this implies that if we choose as $\Omega$ 
the differential form $\omega_{\gamma}$ corresponding to an admissible 
graph, since each $d \phi_{e}$ is obviously closed and the manifolds
$\bar{C}^{+}_{n, \bar{n}}$ are compact by construction, the following holds:
\begin{equation} \label{stokes}
\int_{ \de \bar{C}^{+}_{n, \bar{n}}}\!\! \omega_{\Gamma} = 
\int_{\bar{C}^{+}_{n, \bar{n}}} \!\! d\, \omega_{\Gamma} = 0.
\end{equation}

We will now expand the l.h.s.\ of \eqref{stokes} to show that 
it gives exactly the coefficient $c_{\Gamma}$ 
occurring in \eqref{form2} for the corresponding 
admissible graph.\\

First of all, we want to give an explicit description 
of the manifold $\de \bar{C}^{+}_{n, \bar{n}}$ on which 
the integration is performed. Since the weights $w_{\Gamma}$ 
involved in \eqref{form} are set to vanish identically if the 
degree of the differential form does not match the dimension
of the space on which we integrate, we can restrict our attention
to codimension $1$ strata of $\de \bar{C}^{+}_{n, \bar{n}}$, which have
the required dimension $2 n + \bar{n} - 3$ equal to the number 
of edges and thus of the $1$\ndash forms $d \phi_{e}$.

In an intuitive description of the configuration space 
$\bar{C}_{n,  \bar{n}}$, the boundary represents the degenerate
configurations in which some of the $n+\bar{n}$ points 
``collapse together''. 
The codimension $1$ strata of the boundary can thus be classified as follows:
\begin{itemize}
\item strata of type {\sf S1}, in which $i \geq 2$ points
in the upper half plane $\calH^+$ collapse together
to a point still lying above the real line. Points in
such a stratum can be locally described by the product
\begin{equation} \label{def-S1}
C_{i} \times C_{n-i+1, \bar{n}}.
\end{equation}
where the first term stand for the relative position of the 
collapsing points as viewed ``through a magnifying glass'' 
and the second is the space of the remaining points plus a single
point toward which the first $i$ collapse.

\item strata of type {\sf S2}, in which $i > 0$ points in $\calH^+$
and $j > 0$ points in $\bbR$ with $2 i + j \geq 2$ collapse 
to a single point on the real line. The limit configuration
is given in this case by
\begin{equation} \label{def-S2}
C_{i,j} \times C_{n-i, \bar{n}-j +1}.
\end{equation}
\end{itemize}

\begin{figure}[ht] 
\begin{center}
\resizebox{10 cm}{!}{\includegraphics{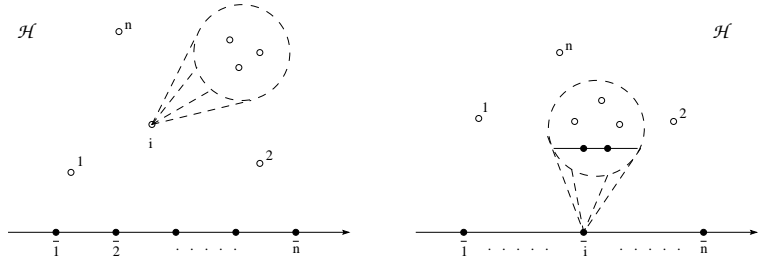}}
\caption{Looking at codimension $1$ strata 
``through a magnifying glass''.} \label{strata}
\end{center}
\end{figure}

These strata have a pictorial representation in Figure \ref{strata}.
In both cases the integral of $\omega_{\Gamma}$ over the stratum can
be split into a product of two integrals of the form 
\eqref{def-weight}: the product of those $d \phi_e$ for which the edge
$e$ connects two collapsing points is integrated over the first
component in the decomposition of the stratum given by 
\eqref{def-S1} resp. \eqref{def-S2}, while the remaining $1$\ndash
forms are integrated over the second.

According to this description, we can split the integral
in the l.h.s.\ of \eqref{stokes} into a sum over different
terms coming from strata of type {\sf S1} and {\sf S2}.
Now we are going to list all the possible configurations
leading to such strata to show that most of these terms vanish and
that the only remaining terms are exactly those 
required to give rise to \eqref{form}. We will not check 
directly that the signs we get by the integration match
with those in \eqref{form}, since we did not give explicitly 
the orientation of the configuration spaces and of their boundaries, 
but we refer once again the reader to the only
paper completely devoted to the careful computation of 
all signs involved in Kontsevich's construction \cite{AMM}.

Among the strata of type {\sf S1}, we distinguish two subcases, according
to the number $i$ of vertices collapsing. Since the integrals 
are set to vanish if the degree of the form does not match the
dimension of the domain, a simple dimensional argument shows 
that the only contributions come from those graphs $\Gamma$
whose subgraph $\Gamma_1$ spanned by the collapsing vertices contains 
exactly $2 i -3 $ edges.

If $i=2$ there is only an edge $e$ involved and in 
the first integral coming from the decomposition 
\eqref{def-S1} the differential of the angle function is integrated
over $C_2 \cong S^1$ and we get (up to a sign) a factor $2 \pi$ which 
cancels the coefficient in \eqref{def-weight}. The  remaining 
integral represents the weight of the corresponding quotient graph $\Gamma_2$
obtained from the original graph after the contraction of $e$: 
to the vertex $j$ of type I resulting from this contraction 
is now associated the $j$\ndash composition of the two 
multivector fields that were associated to the endpoints of
$e$. Therefore, summing over all graphs and all strata of this subtype 
we get the r.h.s.\ of the desired equation \eqref{form}.

If $i\geq 3$, the integral corresponding to this stratum 
involves the product of $2 i -3$ angle forms over $C_{i}$ and
vanishes according to the following Lemma, which contains 
the most technical result among Kontsevich's ``vanishing theorems''.

\begin{figure}[ht] 
\begin{center}
\resizebox{12 cm}{!}{\includegraphics{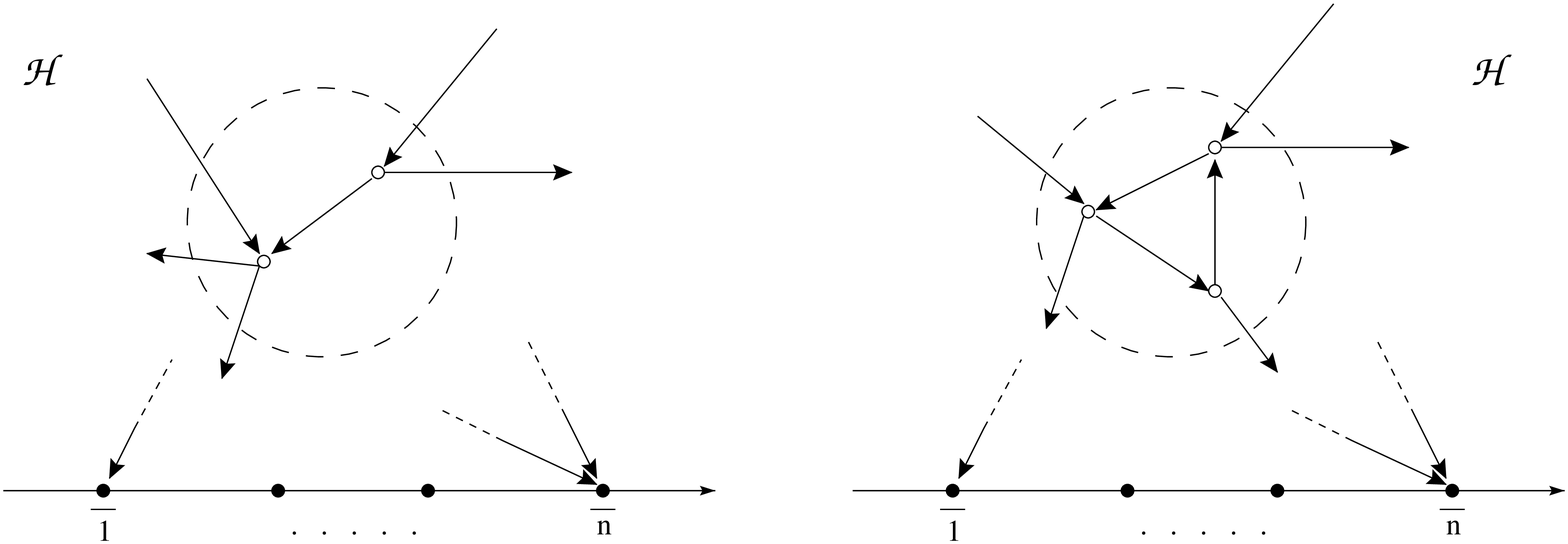}}
\caption{Example of a non vanishing and of a 
vanishing term.} \label{case-S1}
\end{center}
\end{figure}

The two possible situations are exemplified in Figure \ref{case-S1}.

\begin{Lem}
The integral over the configuration space $C_n$ of $n \geq 3$
points in the upper half plane of any $ 2 n -3$ ($=\dim C_n$) 
angle forms $d \phi_{e_i}$ with $i=1, \ldots n$ vanishes 
for $n \geq 3$ 
\end{Lem}
\begin{proof}
The first step consists in restricting the integration
to an even number of angle forms. This is achieved by identifying
the configuration space $C_n$ with the subset of $\calH^n$
where one of the endpoints of $e_1$ is set to be
the origin and the second is 
bounded to lie on the unit circle (this particular configuration
can always be achieved with the help of the action of the
Lie group $G'$). The integral decomposes then into a product 
of $d \phi_{e_1}$ integrated over $S^1$ and the remaining 
$2 n-4 =: 2 N$ forms integrated over the resulting 
complex manifold $U$ given by the isomorphism 
$C_n \cong S^1 \times U$. 
The claim is then a consequence of the following 
chain of equalities:
\begin{equation} \label{trick}
\begin{aligned}
\int_{U} \bigwedge_{j=1}^{2N} d \arg(f_j) =& 
\int_{U} \bigwedge_{j=1}^{2 N} d \log |f_j| = 
\int_{\overline{U}} \calI \Big( d \big( \log |f_1| 
\bigwedge_{j=2}^{2 N} d \log |z_j|\big)\Big)=\\
=& \int_{\overline{U}} d \calI \Big( \big( \log |f_1| 
\bigwedge_{j=2}^{2 N} d \log |z_j|\big)\Big) = 0
\end{aligned}
\end{equation}
where we gave an expression for the angle function
$\phi_{e_j}$ in terms of the argument of the (holomorphic)
function $f_j$ (which is nothing but the difference of 
the coordinates of the endpoints of $e_j$).

The first equality is what Kontsevich calls 
a ``trick using logarithms'' and follows from the decompositions
\[
d \arg(f_j) = \frac{1}{2 i}\, 
\big( d \log (f_j) - d \log (\overline{f}_j) \big)
\]
and
\[
d \log|f_j| = \frac{1}{2}\, 
\big( d \log (f_j) + d \log (\overline{f}_j) \big).
\]
The product of $2 N$ such expressions is thus a 
linear combination of products of $k$ holomorphic
and $2 N- k$ anti-holomorphic forms.
A basic result in complex analysis ensures that,
upon integration over the complex manifold $U$,
the only terms that do not vanish  are those with $k = N$.
It is a straightforward computation to check 
that the non vanishing terms coming from the first
decomposition match with those coming from the second.

In the second equality the integral of the differential 
form is replaced by the integration of a suitable 
$1$\ndash form with values in the space of distributions 
over the compactification $\overline{U}$ of $U$.
A final Lemma in \cite{Ko2} shows that this 
map $\calI$ from standard to distributional
$1$\ndash forms commutes with the differential,
thus proving the last step in \eqref{trick}.
In \cite{Kho}, Khovanskii gave a more elegant proof
of this result in the category of complete complex 
algebraic varieties, deriving the first equality rigorously
on the set of non singular points of $X$ and resolving the 
singularities with the help of a local representation in 
polar coordinates.
\end{proof}

Finally, turning our attention to the strata of type {\sf S2},
the same dimensional argument introduced for the
previous case restricts the possible non vanishing 
terms to the condition that the subgraph $\Gamma_1$ 
spanned by the $i + j$ collapsing vertices (resp. of the first 
and of the second type) contains exactly $2 i + j -2$ edges.

With the same definition as before for the quotient graph
$\Gamma_2$ obtained by contracting $\Gamma_1$, we 
claim that the only non vanishing contributions
come from those graphs for which both graphs
obtained from a given $\Gamma$
are admissible. In this case the weight $w_{\Gamma}$
will decompose into the product $w_{\Gamma_1} \cdot w_{\Gamma_2}$
which in general, by the conditions on the number of edges of $\Gamma$
and $\Gamma_1$,  does not vanish.

\begin{figure}[ht] 
\begin{center}
\resizebox{12 cm}{!}{\includegraphics{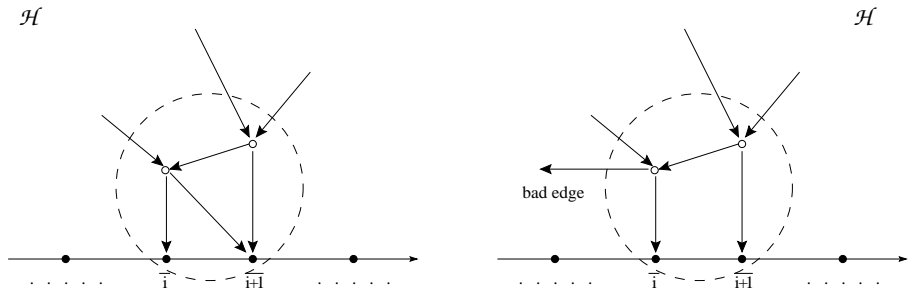}}
\caption{Example of a collapse leading to an admissible quotient
graph and of a collapse correspondidng to a  vanishing term 
because of a bad edge.} \label{case-S2}
\end{center}
\end{figure}

Since all other properties 
required by Definition \ref{def-adm-graphs}
are inherited from $\Gamma$, we have only to check that we do 
not get ``bad edges'' by contraction. The only such possibility
is depicted in the graph on the right in Figure \ref{case-S2}
and occurs when $\Gamma_2$ contains an edge which starts from 
a vertex of the second type: in this case 
the corresponding integral vanishes because it contains
the differential of an angle function evaluated on the pair 
$(z_1, z_2)$, where the first point is constrained to lie
on the real line and such a function vanishes for every $z_2$
because the angle is measured w.r.t.\ the Poincar\'e metric 
(as it can be inferred intuitively from Figure \ref{angle}).

The only non vanishing terms thus correspond 
to the case when we plug the differential operator 
corresponding to the subgraph $\Gamma_1$
as $k$\ndash th argument of the one corresponding
to $\Gamma_2$, where $k$ is the vertex of the second type 
emerging from the collapse. Summing over all such possibilities
and having checked (up to a sign as usual) 
that we get the right weights, it should be clear that the 
contribution due to the strata of type {\sf S2} 
accounts for the l.h.s.\ of \eqref{form}.

In conclusion, we have proved that the morphism $U$ 
is an \LL morphism and since its first coefficient $U_1$
coincides with the map $U_1^{(0)}$ given in Section \ref{hkr-U1}
it is also a quasi-isomorphism and thus determines uniquely 
a star product given by \eqref{U-pi} for any given 
bivector field $\pi$ on $\bbR^d$.

\section{From local to global deformation quantization} \label{globalization}

The content of this last section is based mainly on the work 
of Cattaneo, Felder and Tomassini \cite{CFT1} (see also \cite{CFT2} and \cite{CF2}),
who gave a direct construction of the quantization of
a general Poisson manifold. 

The Kontsevich formula, in fact,
gives a quantization only for the case $M = \bbR^d$ for any Poisson
bivector field $\pi$ and can thus be adopted in the general case to
give only a local expression of the star product. 

The globalization Kontsevich sketched in \cite{Ko2} was carried through
in \cite{Ko3} by abstract arguments, extending the formality theorem to the 
general case.

The works of Cattaneo, Felder and Tomassini instead give 
an explicit recipe to define the star product globally, 
in a similar way to what Fedosov has done in the symplectic
category \cite{Fed}. Also in their approach, the main tool
is a flat connection $\Dbar$ on a vector bundle over $M$ such that
the algebra of the horizontal sections w.r.t.\ to $\Dbar$
is a quantization of the Poisson algebra of the manifold.

We give now an outline of the construction,
addressing the reader to \cite{CFT1} for details and proofs.

In the first step, we introduce the vector bundle $E_0 \to M$
of infinite jets of functions together with the canonical
flat connection $D_0$. The fiber $E_0^x$ over
$x \in M$ is naturally a commutative algebra and inherits the 
Poisson structure induced fiberwise by the Poisson structure on
$\Cinf{M}$. The canonical map which associates to any globally 
defined function its infinite jet at each point $x$ is a Poisson
isomorphism onto the Poisson algebra of horizontal sections 
of $E_0$ w.r.t.\ $D_0$.

As the star product yields a deformation of the pointwise
product on $\Cinf{M}$, we need also a ``quantum version'' of 
the vector bundle and of the flat connection in order to 
find an analogous isomorphism. The vector bundle
$E \to M$ is defined in terms of a section $\aff{\phi}$
of the fiber bundle $\aff{M} \to M$, where $\aff{M}$ is the quotient
of the manifold $\coor{M}$ of jets of coordinates systems 
on $M$ by the action of the group $\GL(d, \bbR)$ of linear diffeomorphisms,
namely $E := (\aff{\phi})^* \widetilde{E}$ where
$\widetilde{E}$ is the bundle of $\bbR\eps$-modules
\[
\coor{M} \; \times_{\GL(d, \bbR)}\; \bbR[[y^1, \ldots , y^d]]\eps \to
\aff{M}.
\]
Since the section $\aff{\phi}$ can be realized explicitly
by a collection of infinite jets at $0$ of maps 
$\phi_x \colon \bbR^d \to M$ such that $\phi_x(0) = x$
for every $x \in M$ (defined modulo the action of $GL(d, \bbR)$),
we can suppose for simplicity that we have fixed 
a representative $\phi_x$ of the equivalence class for
each open set of a given covering, thus realizing a trivialization
of the bundle $E$. Therefore, from now on we will
identify $E$ with the trivial bundle with fiber 
$ \bbR[[y^1, \ldots , Y^d]]\eps$; in this way $E$ realizes 
the desired quantization, since it is isomorphic (as a bundle
of $\bbR\eps$-modules) to the bundle $E_0\eps$ whose elements
are formal power series with infinite jets of functions as coefficients.

In order to define the star product and the connection
on $E$, we have first to introduce some new objects 
whose existence and properties are byproducts of the formality
theorem.
Given a Poisson bivector field $\pi$ and two vector fields
$\xi$ and $\eta$ on $\bbR^d$, we define:
\begin{equation} \label{def-P-A-F}
\begin{aligned}
P(\pi) &:= \sum_{k=0}^{\infty} \frac{\ee^k}{k!}\; U_k (\pi, \dots, \pi),\\
A(\xi, \pi) &:= \sum_{k=0}^{\infty} \frac{\ee^k}{k!}\; U_{k+1} (\xi, \pi, \dots, \pi),\\
F(\xi,\eta,\pi) &:= \sum_{k=0}^{\infty} \frac{\ee^k}{k!}\; U_{k+2} (\xi,\eta,\pi, \dots, \pi).
\end{aligned}
\end{equation}

A straightforward computation of the degree of the multidifferential
operators on the r.h.s.\ of \eqref{def-P-A-F} shows that $P(\pi)$ is
a (formal) bidifferential operator, $A(\xi, \pi)$ a differential
operator and $F(\xi,\eta,\pi)$ a function. Indeed $P(\pi)$ is 
nothing but the star product associated to $\pi$ as introduced at the
end of Section \ref{digression}.

More precisely, $P$, $A$ and $F$ are elements of degree
resp. $0$, $1$ and $2$ of the Lie algebra cohomology complex
of (formal) vector fields with values in the space of local polynomial 
maps, i.e.\ multidifferential operators depending polynomially on $\pi$: 
an element of degree $k$ of this complex is a map 
that sends $\xi_1 \wedge \cdots \wedge \xi_k$ to 
a multidifferential operator $S(\xi_1,\ldots , \xi_k, \pi)$
(we refer the reader to \cite{CFT1} for details).
The differential $\delta$ on this complex is then 
defined by
\begin{equation}
\begin{aligned}
\delta\, S(\xi_1,\ldots , \xi_{k+1}, \pi) := 
&\sum_{i=1}^{k+1} (-)^i \frac{d}{d t} \Big|_{t=0} 
S(\xi_1,\ldots, \hat{\xi}_i, \ldots, \xi_{k+1},(\Phi_\xi^t)_*\, \pi)\; + \\
+ & \sum_{i<j}  (-)^{i+j} \;
S(\Lie{\xi_i}{\xi_j},\xi_1,\ldots, \hat{\xi}_i, \ldots, \hat{\xi}_j, x\ldots, \xi_{k+1}, \pi)
\end{aligned}
\end{equation}
where a caret denotes as usual the omission of the corresponding argument
and $\Phi_\xi^t$ is the flow of the vector field $\xi$.

As the associativity condition on the star product, which can now be written 
in the form $P \circ (P \otimes \id - \id \otimes P) = 0$, follows
from the formality theorem, the following equations are a corollary 
of the same result and can be proved with analogous computations:
\begin{itemize} 
\item 
$P(\pi) \circ (A(\xi,\pi) \otimes \id + \id \otimes A(\xi,\pi))
= A(\xi, \pi) \circ P(\pi) + \delta P(\xi,\pi)  $
\item 
$P(\pi) \circ ( F(\xi,\eta,\pi) \otimes \id - \id \otimes F(\xi,\eta,\pi))= $
\begin{flushright}
$= A(\xi,\pi) \circ A(\eta,\pi) - A(\eta,\pi) \circ A(\xi,\pi) + \delta A(\xi,\eta,\pi)$
\end{flushright}
\item $-A(\xi,\pi) \circ F(\eta,\zeta,\pi) -A(\eta,\pi) \circ
  F(\zeta,\xi,\pi) -A(\zeta,\pi) \circ F(\xi,\eta,\pi)
= \delta F(\xi,\eta,\zeta,\pi)$
\end{itemize}
\begin{equation} \label{P-A-F-rel}\, \end{equation}
The first of these equations describes the fact that under 
the coordinate transformation induced by $\xi$ the star product
$P(\pi)$ is changed to an equivalent one up to higher order terms.
The last two equations will be used
in the construction of the connection and its curvature, since 
they represent an analogous of the defining relations between 
a connection $1$\ndash form $A$ and its curvature $F_A$.

Upon explicit computation of the configuration
space integrals involved in the definition of the 
Taylor coefficients $U_k$, we can also give the lowest order terms
in the expansion of $P$, $A$ and $F$ and their action on functions:
\begin{enumerate}
\item[(i)] $P(\pi)(f\otimes g)= f\, g + \ee\, \pi(df,dg) + O(\ee^2)$;
\item[(ii)] 
$A(\xi,\pi)= \xi+O(\epsilon)$, where we identify $\xi$
with a first order differential operator on the r.h.s.;
\item[(iii)] $A(\xi,\pi)=\xi$, if $\xi$ is a linear
vector field;
\item[(iv)]
$F(\xi,\eta,\alpha)=O(\ee)$;
\item[(v)] $P(\pi)(1\otimes f)=P(\pi)(f\otimes1)=f$;
\item[(vi)] $A(\xi,\pi)1=0$.
\end{enumerate}

Equations $i)$ and $v)$ where already introduced 
in Definition \ref{def-star-prod} as two of the defining 
conditions of a star product, while the ones involving
$A$ are used to construct a connection $D$ on sections
of $E$.

A section $f \in \Gamma(E)$ is given locally by a map
$x \to f_x$ where for every $y$, $f_x(y)$ is a formal power 
series whose coefficients are infinite jets. On the space of
such sections we can introduce a deformed product $\star$ 
which will give us the desired star product on $\Cinf{M}$
once we identify horizontal sections with ordinary functions.
Denoting analogously by $\pi_x$ the push-forward by $\phi_x^{-1}$
of the Poisson bivector $\pi$ on $\bbR^d$, we can define
the deformed product through the formal bidifferential 
operator $P(\pi_x)$ in the same way as $P(\pi)$ represents
the usual star product:
\[
(f \star g )_x (y) := f_x(y)\;g_x(y) + 
\ee\; \pi_x^{i j}(y) \frac{\de f_x(y)}{\de y^i} \frac{\de g_x(y)}{\de y^j} 
+ O(\ee^2).
\]
We can define the connection $D$ on $\Gamma(E)$
by 
\[
(D \, f)_x = d_x f + A_x^{M} f
\] 
where $d_x f$ is the de Rham differential of $f$ regarded 
as a function with values in $\bbR[[y^1, \ldots, y^d]]\eps$ 
and the formal connection $1$\ndash form 
is specified by its action on a tangent vector $\xi$ by
\[
A_x^{M}(\xi) = A(\hat{\xi}_x, \pi_x)
\]
where $A$ is the operator defined in \eqref{def-P-A-F}
evaluated on the multivector fields $\xi$ and $\pi$
expressed in the local coordinates system given by $\phi_x$.

The important point is that since the coefficients $U_k$
of the formality map that appear in the definition 
of $P$ and $A$ are polynomial in the derivatives of
the coordinate of the arguments $\xi$ and $\pi$, 
all results holding for $P(\pi)$ and $A(\xi,\pi)$
are inherited by their formal counterparts. In particular
equalities $i)$ and $v)$ above (together with the formality 
theorem from which they are derived) ensure that $\star$
is an associative deformation of the pointwise 
product on sections and equalities $ii)$ and $iii)$
can be used to prove that $D$ is indeed independent
of the choice of $\phi$ and therefore
induces a global connection on $E$.

We can finally extend $D$ and $\star$ by
the (graded) Leibniz rule to the whole complex of
formal differential forms 
$\Omega^{\bullet}(E) = \Omega{M} \otimes_{\Cinf{M}} \Gamma(E)$
and use \eqref{P-A-F-rel} to verify the following
\begin{Lem}
Let $F^{M}$ be the $E$\ndash valued $2$\ndash form
given by $x \to F_x^{M}$ 
where $F_x^{M}(\xi, \eta) = F(\hat{\xi}_x, \hat{\eta}_x, \pi_x)$
for any pair of vector fields $\xi$, $\eta$.
Then $F^{M}$ represent the curvature of $D$ and 
the two are related to each other and to the star product
by the usual identities:
\begin{enumerate}
\item[a)] $D( f \star g) =D(f) \star g + f \star D(g)$;
\item[b)] $D^2(\cdot) = \starLie{F^{M}}{\cdot}$;
\item[c)] $D\, F^{M} = 0$
\end{enumerate}
\end{Lem}
\begin{proof}
The identities follow directly from the relations
\eqref{P-A-F-rel}, in which the star commutator 
$\starLie{f}{g}= f \star g - g \star f$
is already implicitly defined,
once we identify the complex of formal multivector fields
endowed with the differential $\delta$ with the complex
of formal multidifferential operators with the de Rham 
differential. The map that realizes this isomorphism
is explicitly defined in \cite{CFT1}.
\end{proof}

A connection $D$ satisfying the above relations
on a bundle $E$ of associative algebras is called 
a {\sf Fedosov connection} with {\sf Weyl curvature}
$F$: it is the kind of connection Fedosov introduced
to give a global construction in the symplectic
case. Following Fedosov, the last step to 
the required globalization is to deform $D$ into
a new connection $\Dbar$ which enjoys the same properties
and moreover has zero Weyl curvature, so that we
can define the complex $H^k(E, \Dbar)$ and in
particular the (sub)algebra of horizontal sections
$H^0(E, \Dbar)$.

The construction of $\Dbar$ relies on the following Lemmata.

\begin{Lem}
Let $D$ be a Fedosov connection on $E$ with Weyl curvature 
$F$ and $\gamma$ an $E$\ndash valued $1$\ndash form, then
\[
\Dbar := D + \starLie{\gamma}{\cdot}
\]
is also a  Fedosov connection whose Weyl curvature is
$\overline{F} = F + D\, \gamma + \gamma \star \gamma$.
\end{Lem} 
\begin{proof}
For any given section $f$, a direct computation shows
\[
\begin{aligned}
\Dbar^2 f &= \starLie{F}{f} + D \starLie{\gamma}{f} 
+ \starLie{\gamma}{D f} + \starLie{\gamma}{\starLie{\gamma}{f}} =\\
&= \starLie{F}{f} + \starLie{D \gamma}{f} + \starLie{\gamma}{\starLie{\gamma}{f}} =
\starLie{F + D\, \gamma + \frac{1}{2} \starLie{\gamma}{\gamma}}{f}
\end{aligned}
\]
where the last equality follows from the Jacobi 
identity for the star commutator, since every associative
product induces a Lie bracket given by the commutator.

Applying $\Dbar$ on the new curvature, we can check explicitly that
\[
\begin{aligned}
\Dbar \Big(F + D\, \gamma + \frac{1}{2} \starLie{\gamma}{\gamma}\Big)
&= D^2 \gamma + \frac{1}{2} \starLie{D \gamma}{\gamma} -
\frac{1}{2} \starLie{\gamma}{D \gamma} + \starLie{\gamma}{F+ D\gamma}=\\
&= \starLie{F}{\gamma} + \starLie{\gamma}{F} = 0
\end{aligned}
\]
where we made use again of the (graded) Jacobi identity and of
the (graded) skew-symmetry of $\starLie{}{}$.
\end{proof}

\begin{Lem}
Let $D$ be a Fedosov connection on a bundle $E = E_0\eps$ and $F$
its Weyl curvature and let 
\[
D = D_0 + \ee D_1 + \cdots 
\qquad \text{and} \qquad 
F = F_0 + \ee F_1 + \cdots
\]
be their expansions as formal power series.
If $F_0 =0 $ and the second cohomology of $E_0$ w.r.t.\ $D_0$
is trivial, there exist a $1$\ndash form $\gamma$ such that 
$\Dbar$ has zero Weyl curvature.
\end{Lem} 
\begin{proof}
By the previous Lemma, the claim is equivalent to 
the existence of a solution to the equation
\[
\overline{F} = F + D\, \gamma + \frac{1}{2} \starLie{\gamma}{\gamma} = 0.
\]
A solution can be explicitly constructed by induction on 
the order in $\ee$, starting from $\gamma_0 = 0$
and assuming that $\gamma^{(k)}$ is a solution $\mod \ee^{k+1}$.
We can thus add to 
$\overline{F}^k = F + D\, \gamma^{(k)} + 
\frac{1}{2} \starLie{\gamma^{(k)}}{\gamma^{(k)}}$
the next term $\ee^{k+1} D_0\, \gamma_{k+1}$
to get $\overline{F}^{(k+1)}$ modulo higher terms.
From $D \overline{F}^{(k)} + 
\starLie{\gamma^{(k)}}{\overline{F}^{(k)}} =0$ 
and the induction hypothesis $\overline{F}^{(k)} =0$ 
$\mod  \ee^{k+1}$ we get $D_0\, \overline{F}^{(k)} = 0$.
Since now $H^2(E_0, D_0) = 0$, we can invert $D_0$
to define $\gamma_{k+1}$ in terms of the lower order terms
$\overline{F}^{(k)}$ in such a way that $\overline{F}^{(k+1)} = 0$
is satisfied $\mod \ee^{k+2}$, thus completing the induction step.
\end{proof}

Since in our case $D$ is a deformation 
of the natural flat connection $D_0$ on sections of
the bundle of infinite jets, the hypothesis of the previous 
Lemma are satisfied and we can actually find 
a flat connection $\Dbar$ which is still a good 
deformation of $D_0$.

A last technical Lemma gives us an isomorphism
between the algebra of the horizontal sections
$H^0(E, \Dbar)$ and its non-deformed counterpart
$H^0(E_0, D_0)$, which in turn is isomorphic 
to the Poisson algebra $\Cinf{M}$: this concludes 
the globalization procedure.

Only recently, Dolgushev \cite{Do} gave a new proof
of Kontsevich's formality theorem for a general 
manifold. The main difference in this approach is 
that it is based on the use of covariant tensors unlike 
Kontsevich's original proof, which is based on $\infty$-jets 
of multidifferential operators and multivector fields
and is therefore intrinsically local. 
In particular, he gave a solution of the 
deformation quantization problem for an arbitrary Poisson 
orbifold.

%%%%%%%%%%% BIBLIOGRAPHY %%%%%%%%%%%%%

\end{document}